\documentclass[fleqn,10pt]{wlscirep}
\usepackage[utf8]{inputenc}
\usepackage{algorithm2e}
\usepackage{stix}
\usepackage{multirow}
\usepackage[symbol]{footmisc}
\usepackage{comment}

\usepackage{epsfig}
\usepackage{soul}

\usepackage{tikz}
\usepackage{tkz-fct}
\usetikzlibrary{arrows,shapes,trees}
\def\thickhline{\noalign{\hrule height.9pt}}
\RestyleAlgo{ruled}
\usepackage{titletoc}
\SetLabelAlign{Center}{\hfil#1\hfil}
\SetLabelAlign{CenterWithParen}{\hfil(\makebox[1.0em]{#1})\hfil}

\definecolor{irishgreen}{RGB}{0,154,73}
\newcommand{\af}{\color{irishgreen}}

\usepackage[T1]{fontenc}
\title{Detecting seizure onset and offset times using human intelligence:\\ A critical-transitions-based approach}

\author[1,2,*]{Andrew Flynn}
\author[3]{Cian McCafferty}
\author[4,5,6]{Klaus Lehnertz}
\author[7]{Fran\c{c}ois David}
\author[8]{Vincenzo Crunelli}
\author[2,9]{Gordon Lightbody}
\author[1]{Sebastian Wieczorek}
\affil[1]{School of Mathematical Sciences, University College Cork, T12 XF62 Cork, Ireland.}
\affil[2]{INFANT Research Centre, University College Cork, T12 DC4A Cork, Ireland.}
\affil[3]{Department of Anatomy \& Neuroscience, University College Cork, Cork, Ireland.}
\affil[4]{Department of Epileptology, University of Bonn Medical Centre, 53127 Bonn, Germany.}
\affil[5]{Helmholtz-Institute for Radiation and Nuclear Physics, University of Bonn, 53115 Bonn, Germany.}
\affil[6]{Interdisciplinary Center for Complex Systems, University of Bonn, 53175 Bonn, Germany.}
\affil[7]{Center for Interdisciplinary Research in Biology, Coll{\`e}ge de France, 75005 Paris, France}
\affil[8]{Department of Pharmacology and Neuroscience, Faculty of Medicine, University of Lisbon, Lisbon, Portugal}
\affil[9]{Department of Electrical and Electronic Engineering, University College Cork, T12 YF78 Cork, Ireland.}

\affil[*]{andrewflynn@ucc.ie}

\begin{abstract}

Most existing seizure detection algorithms require extensive pre-processing of the data and rely on heuristic or currently unexplainable machine learning approaches. These approaches often struggle with balancing detection sensitivity and specificity in the presence of variable seizure morphologies, interictal epileptiform discharges, and artefacts.
Here, we consider an alternative approach: our seizure detection algorithm, which is based on the concept of critical transitions and overcomes the aforementioned limitations. 
Specifically, we perform a receiver-operating-characteristic analysis to quantify the performance of our algorithm in terms of its agreement with expert annotations of seizure onset and offset times in the voltage recordings of seizure activity in epileptic rodents with different seizure morphologies.
We demonstrate how performance depends on algorithm parameters and varies across different rodent recording sessions.
We determine the optimal set of algorithm parameters for each recording session, with near expert-level performance achieved in most cases.
Finally, we derive a single general set of algorithm parameters applicable across all recording sessions.
The algorithm maintains its high performance in this general setting, demonstrating its versatility, robustness across varying seizure morphologies, and potential to complement machine learning algorithms.

\end{abstract}
\begin{document}

\flushbottom
\maketitle

\thispagestyle{empty}

\newpage

\section*{Introduction}

Seizure detection algorithms have become a valuable tool in the diagnosis, monitoring, and management of epilepsy, a condition which effects $50$ million people worldwide~\cite{WHO_Seizures}. 
However, despite decades of research and improvements in obtaining voltage recordings of seizure activity from the brain, significant challenges remain in developing robust and accurate seizure detection algorithms. 
Factors that add to this challenge include the susceptibility of voltage recordings to artefacts~\cite{islam2016artefactdetection}, the presence of interictal epileptiform discharges~\cite{decurtis2012interictal,niedermeyerdasilva2005_book}, and the variability of seizure morphology across individuals, recording sites, and recording sessions. 
As a result, many algorithms struggle to balance sensitivity and specificity - correctly identifying seizure and non-seizure intervals - in clinically realistic settings, yet major progress has been made thanks to modern machine learning methods~\cite{Temko2011_Performance,Mathieson2016_Performance}. 
However, the price to pay for using certain machine learning techniques is generally their black-box nature and tendency to confabulate (hallucinate); they lack explainability and generate false, but sometimes plausible, predictions/information\cite{Flynn25_confab}. 
Addressing the above challenges is essential for advancing personalised epilepsy care and realising the full potential of neurotechnological interventions, such as cortical stimulation~\cite{sun2008cortstim,morrell2011cortstim}, which require real-time information on the state of the brain.

In this paper, we conduct a more detailed analysis on the performance of a seizure detection algorithm introduced in our previous work, Flynn \textit{et al.}~\cite{flynn2025classifying}. 
To provide a different viewpoint to and complement recent trends in using artificial intelligence to develop seizure detection algorithms, we use an ancient technology known as `human intelligence'.
To develop this algorithm we first acknowledged that the epileptic seizures we analysed are states of high-amplitude, synchronous electrical activity in the brain with defined onset and offset\cite{beniczky2025_LeagueAgainstEpil_def}. Further, while seizures have a wide range of associated symptoms, seizures constitute a sudden and large change in brain state and have a commensurate impact on the life of a person with epilepsy\cite{kaye2025_seizimpact}. 
From a mathematical perspective, a sudden and large change in the state of a complex system corresponds with Ashwin~\textit{et al.}'s~\cite{ashwinwieczorek2017CTdef} definition of a critical transition (CT).
The mathematical theory of CTs has been instrumental in describing and foreseeing sudden and large changes in climate and ecological systems~\cite{lenton2008tippingpolicymaking,lenton2020tippingpositivechange} and more recently in shaping global policymaking~\cite{armstrong2022exceeding,lenton2025global}. 
Inspired by its impact on the environmental sciences, we consider the brain as a complex system and apply the theory of CTs to detect seizure activity in the brain. Specifically, we model seizure onset and offset as CTs between a non-seizure state (NS state) and seizure state (S state), and develop an algorithm to detect CTs between these states in voltage recordings of brain activity. Importantly, our algorithm is based entirely on mathematical definitions of these states and CTs between them. These definitions are informed from observations of seizure activity in local field potential (LFP) recordings but generalisable to electroencephalography (EEG) recordings. Our algorithm detects CTs through monitoring only a single feature of the data, the voltage recording at a given time, i.e., we simply use the data we are given as input to our algorithm. 

In our previous work, Flynn \textit{et al.}~\cite{flynn2025classifying}, we introduced our algorithm as part of a framework for classifying seizure onset in terms of different types of CT. We applied this framework to voltage recordings of seizure activity in \textit{Genetic Absence Epilepsy Rats from Strasbourg} (GAERS), specifically measurements of the LFP in mV. 
In the present paper, we conduct a more detailed analysis on the performance of our algorithm when applied to the same set of voltage recordings; see part \hyperref[M1]{M1} of the Methods section for details on how these voltage recordings are obtained and annotated by an expert. These voltage recordings provide a reasonable challenge for any seizure detection algorithm given there are numerous artefacts and interictal epileptiform discharges present and seizure morphology varies significantly between GAERS. In particular, the change in voltage amplitude from the NS to S state varies between GAERS, with some exhibiting relatively significantly larger changes than others.

In practice, a well-designed seizure detection algorithm will identify seizure and non-seizure time intervals that agree with expert annotations. 
However, one must carefully choose what performance metrics are used to quantify this agreement. 
In this paper we take inspiration from Temko \textit{et al.}\cite{Temko2011_Performance} and Mathieson \textit{et al.}\cite{Mathieson2016_Performance} by using performance metrics based on receiver-operating-characteristic (ROC) curves. 
These metrics enable us to determine the optimal choice of algorithm parameters for each recording session we analyse, with very good performance achieved on average. Based on this we derive a general set of algorithm parameters applicable across all recording sessions. We show that the algorithm maintains much of its accuracy in this general setting, demonstrating its robustness across varying seizure morphologies.

Overall, our results demonstrate the benefits of studying the epileptic brain from a CT perspective. We produce a new approach for accurate and automated seizure detection that is grounded in mathematics and can be adapted depending on definitions of NS and S states. 
Since our approach does not rely on heuristic or currently unexplainable machine learning approaches, it could complement these approaches, improving their overall performance beyond what can be achieved by adjusting the parameters of the machine learning algorithm.

\section*{Results}

\subsection*{Detecting seizure activity in voltage recordings \label{sssec:UsingCTSeizAlgorithm}}

In this subsection we briefly describe the main premise and technical aspects of our algorithm which are needed to discuss our results. 

\subsubsection*{Main premise of our algorithm}
We designed our algorithm through observation of generic characteristics and expert annotations of seizure activity. 
We use Fig.~\ref{fig:RatData_Properties_} from the Methods section as an example where, according to an expert, seizure onset happens at $t = \tau_{1} \approx8.5\text{s}$ and offset at $t = \tau_{2} \approx22\text{s}$. 
From a CT perspective, we observe:
\begin{itemize}
    \item[(O1)\label{O1}] a non-seizure state (NS){\af,} characterised by a small-amplitude and weakly-correlated type of fluctuation around $0~\text{mV}$, 
    \item[(O2)\label{O2}] a seizure state (S){\af,} characterised by a large-amplitude and strongly-correlated type of fluctuation around $0~\text{mV}$, and 
    \item[(O3)\label{O3}] two CTs between the NS and S states, characterised by a sudden change from one type of fluctuation to another. 
\end{itemize}
We consider these two CTs in Fig.~\ref{fig:RatData_Properties_} to be analogous to seizure onset and offset, respectively. This approach is motivated by similarities between the International League Against Epilepsy's definition of a seizure\cite{beniczky2025_LeagueAgainstEpil_def} and Ashwin \textit{et al.}'s definition of a CT\cite{ashwinwieczorek2017CTdef}.

\subsubsection*{Main technical aspects of our algorithm}
Inspired by a `non-ideal relay'~\cite{Pokrovskii12systemswHysteresis}, we define CTs in terms of successive crossings of two voltage thresholds and use a `moving window analysis' to assess how long the brain remains in the new state after successive crossings of the two voltage thresholds. 
By using two thresholds we avoid the shortcomings of a single threshold approach in the case of CTs to the S state; a single threshold can be crossed multiple times within a short time interval given the high frequency of voltage measurements or, naturally, when in the S state given its large-amplitude oscillatory-like nature, both of which could result in false detections.

\begin{figure}[t]
    \centering
    \includegraphics[width=\textwidth]{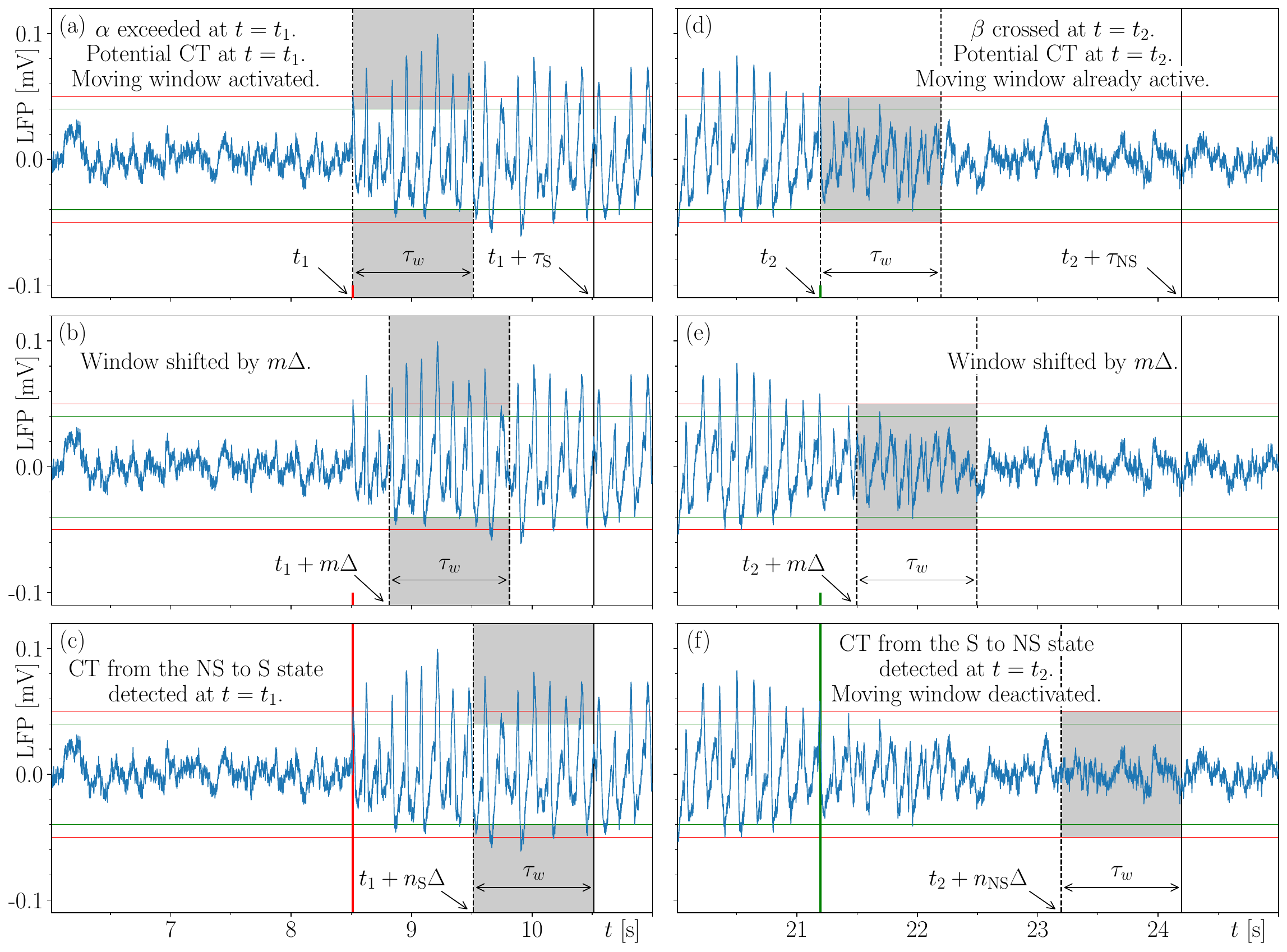}
    \caption{
    Illustrating how the CT detection algorithm works, using the voltage recordings from Fig.~\ref{fig:RatData_Properties_} as an example. (a)-(c) Detecting a CT from the NS to S state at time $t=t_{1}$. (d)-(f) Detecting a CT from the S to NS state at time $t=t_{2}$. The algorithm parameters are chosen as $\alpha = 0.05$, $\beta = 0.04$, $\tau_{\text{S}} = 2$, $\tau_{\text{NS}}=3$, $\tau_{w}=1$, and we use $m=300$ in (b) and (e). The red and green horizontal lines indicate the thresholds of $\alpha$ and $\beta$, respectively.}
    \label{fig:Algorithm_in_use_}
\end{figure}

Briefly, our algorithm detects a CT from the NS to S state at time $t=t_{1}$ if the absolute value of the time series exceeds an upper voltage threshold $\alpha$ and then continues to exceed a lower voltage threshold $\beta$ frequently enough for a period of at least $\tau_{\text{S}}$. Similarly, our algorithm detects a CT from the S to NS state at time $t=t_{2}$ if the absolute value of the time series falls below $\beta$ and then does not exceed $\alpha$ for a period of at least $\tau_{\text{NS}}$. 
A window of length $\tau_{w}$ is moved along the time series in discrete steps of size $\Delta$ to check whether the thresholds are exceeded in each window throughout the specified period. 
Part~\hyperref[M2]{M2} of the Methods section provides precise definitions of when the algorithm detects a CT, specifies the algorithm's technical details, and how the algorithm is designed to mitigate the influence of artefacts and interictal epileptiform discharges. 
In Fig.~\ref{fig:Algorithm_in_use_} we use the actual voltage recording from Fig.~\ref{fig:RatData_Properties_} to illustrate in three simple steps how our algorithm detects a CT in voltage recordings. 
\\
\\
\noindent
The remainder of the present paper is focused on tuning the parameters of our CT detection algorithm to maximise its agreement with the expert's annotations in terms of the corresponding seizure and non-seizure time intervals. 
We denote seizure and non-seizure time intervals according to the expert by $I_{\text{S}}^{(E)}$ and $I_{\text{NS}}^{(E)}$. 
The ends of these intervals are defined by the values of seizure onset times $\tau_{1}$ and seizure offset times $\tau_{2}$. 
Similarly, we use $I_{\text{S}}^{(A)}$ and $I_{\text{NS}}^{(A)}$ to denote the seizure and non-seizure time intervals obtained from the algorithm, by detecting the times of different CTs, namely $t_{1}$ and $t_{2}$; 
see part~\hyperref[M3]{M3} of the Methods section for precise definitions. 
In a given recording session, we denote the number of $I_{\text{S}}^{(E)}$ intervals by $N_{\text{S}}^{(E)}$, and the number of $I_{\text{NS}}^{(E)}$ intervals by $N_{\text{NS}}^{(E)}$, where $N_{\text{NS}}^{(E)} = N_{\text{S}}^{(E)}-1$. We use a similar notation for the algorithm, namely $N_{\text{S}}^{(A)}$ and $N_{\text{NS}}^{(A)}$.

\subsection*{Classification of time intervals, constructing ROC curves, and algorithm performance metrics}

In this paper we take inspiration from Temko \textit{et al.}\cite{Temko2011_Performance} and Mathieson \textit{et al.}\cite{Mathieson2016_Performance} and evaluate the performance of our algorithm using metrics based on \textbf{receiver-operating-characteristic (ROC)} curves, a long standing method from signal detection theory to evaluate sensitivity and specificity, introduced in the 1940s~\cite{fawcett2006introduction}. 
The points on our ROC curves are derived from classifications of $I_{\text{S}}^{(A)}$ and $I_{\text{NS}}^{(A)}$ where the $I_{\text{S}}^{(E)}$ and $I_{\text{NS}}^{(E)}$ are taken as the ground truth. 
In this subsection, we briefly describe the time interval classification procedure, how we construct ROC curves using this procedure, and the performance metrics we obtain from ROC curves.

\subsubsection*{Time interval classification procedure}

In short, a given $I_{\text{S}}^{(A)}$ is classified as \textbf{\textit{true positive} (TP)} if it overlaps with at least one $I_{\text{S}}^{(E)}$. Otherwise $I_{\text{S}}^{(A)}$ is classified as \textbf{\textit{false positive} (FP)}. Similarly, a given $I_{\text{NS}}^{(A)}$ is classified as \textbf{\textit{true negative} (TN)} if it overlaps with at most one $I_{\text{NS}}^{(E)}$. Otherwise $I_{\text{NS}}^{(A)}$ is classified as \textbf{\textit{false negative} (FN)}. 
We introduce two parameters, $\rho_{\text{S}}, \,\rho_{\text{NS}} \in \left[0,1\right]$, to quantify how much overlap is required to classify a given $I_{\text{S}}^{(A)}$ as TP and $I_{\text{NS}}^{(A)}$ as TN.
When $\rho_{\text{S}} = \,\rho_{\text{NS}} = 0$ any overlap between intervals of $I_{\text{S}}^{(A)}$ and $I_{\text{S}}^{(E)}$, and $I_{\text{NS}}^{(A)}$ and $I_{\text{NS}}^{(E)}$, will suffice to classify $I_{\text{S}}^{(A)}$ as a TP or $I_{\text{NS}}^{(A)}$ as a TN, 
whereas when $\rho_{\text{S}} = \rho_{\text{NS}} = 1$ the intervals need to completely overlap, i.e., an $I_{\text{S}}^{(A)}$ needs to contain an $I_{\text{S}}^{(E)}$ or vice-versa, and similarly for $I_{\text{NS}}^{(A)}$ and $I_{\text{NS}}^{(E)}$.
See part~\hyperref[M3]{M3} of the Methods section for a more precise description of our classification procedure.

\begin{figure}[t]
    \centering
    \includegraphics[width=.95\textwidth]{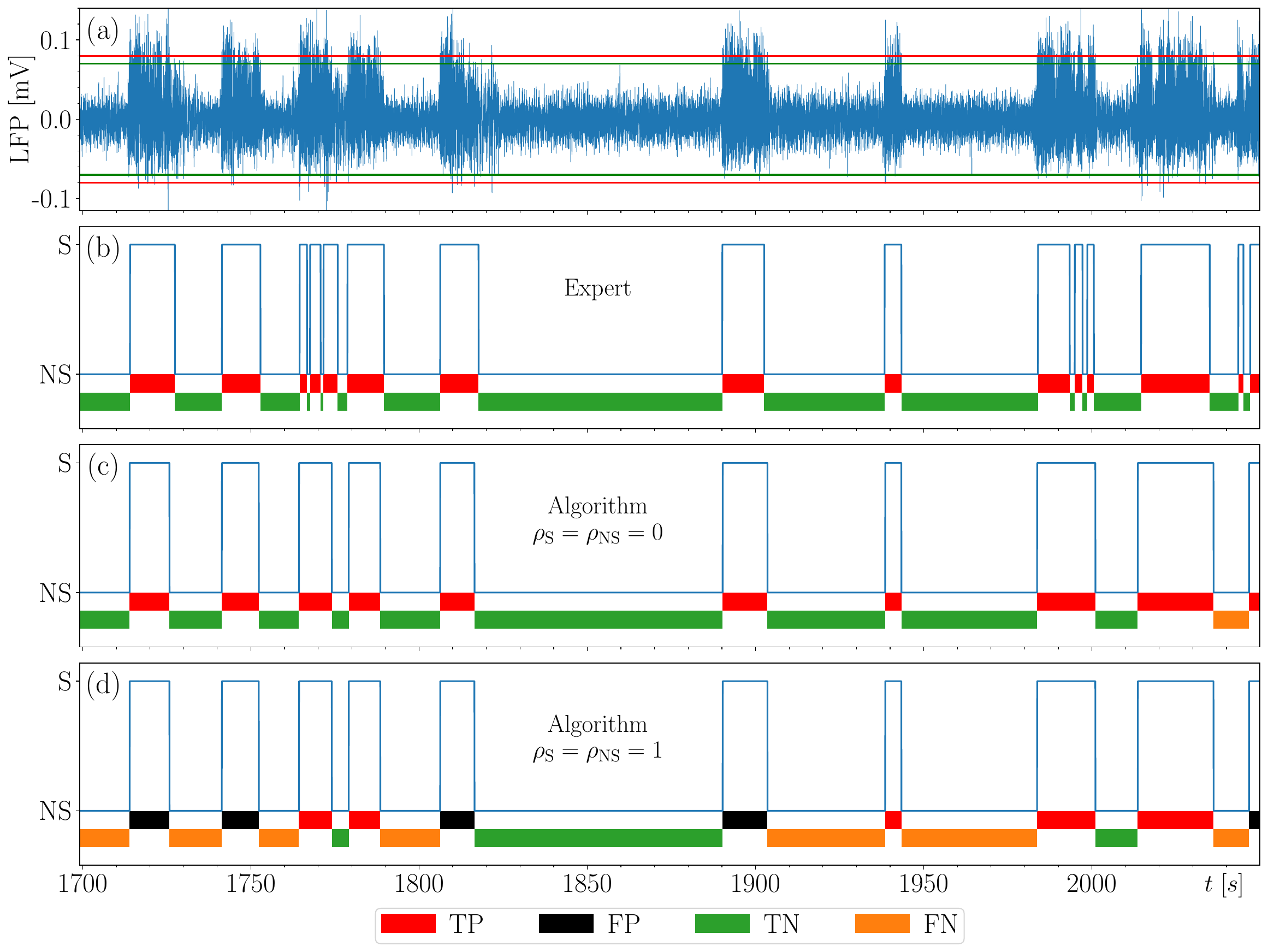}
    \caption{Illustrating how the seizure and non-seizure time interval classification procedure specified in part \hyperref[M3]{M3} of the Methods section works and the effect that $\rho_{\text{S}}$ and $\rho_{\text{NS}}$ have on the classifications. (a) shows a portion of voltage recordings that were annotated by the expert. (b) shows a binary sequence representation of the non-seizure and seizure time intervals according to the expert. (c) and (d) show the corresponding sequence of the non-seizure and seizure time intervals according to the algorithm. The boxes beneath these sequences are coloured green for TN, red for TP, orange for FN, and black for FP. The algorithm's parameters were chosen as $\tau_{\text{NS}}=3$, $\tau_{\text{S}}=2$, $\tau_{w}=1$, with $\alpha=0.08$ and $\beta=0.07$ (indicated by the horizontal red and green lines in (a)), $\rho_{\text{S}} = \rho_{\text{NS}}=0$ in (c), $\rho_{\text{S}} = \rho_{\text{NS}}=1$ in (d).}
    \label{fig:rhoS_effect_ROC_example1}
\end{figure}

Figure~\ref{fig:rhoS_effect_ROC_example1} illustrates how our classification procedure works. 
Specifically, we pick a portion of the voltage recordings from recording session T2M (shown in Fig.~\ref{fig:rhoS_effect_ROC_example1}~(a)) that was annotated by the expert. These annotations are shown in Fig.~\ref{fig:rhoS_effect_ROC_example1}~(b) via a binary sequence representation of the $I_{\text{S}}^{(E)}$ and $I_{\text{NS}}^{(E)}$ with coloured boxes beneath indicating these intervals as the ground truth. 
We use our algorithm to detect CTs between NS and S states in Fig.~\ref{fig:rhoS_effect_ROC_example1}~(a) with algorithm parameters chosen as follows: $\alpha = 0.08$, $\beta = 0.07$, $\tau_{\text{NS}}=3$, $\tau_{\text{S}}=2$, and $\tau_{w}=1$. 
We show the resulting classification of $I_{\text{S}}^{(A)}$ and $I_{\text{NS}}^{(A)}$ via the coloured boxes beneath the binary sequence representation of these time intervals for $\rho_{\text{S}}=\rho_{\text{NS}}=0$ in Fig.~\ref{fig:rhoS_effect_ROC_example1}~(c) and $\rho_{\text{S}}=\rho_{\text{NS}}=1$ in Fig.~\ref{fig:rhoS_effect_ROC_example1}~(d). By comparing Fig.~\ref{fig:rhoS_effect_ROC_example1}~(c) to (d) we see the main effect that $\rho_{\text{S}}$ and $\rho_{\text{NS}}$ have on the classifications - more time intervals are classified as TPs and TNs when $\rho_{\text{S}}=\rho_{\text{NS}}=0$ as opposed to $\rho_{\text{S}} = \rho_{\text{NS}}=1$.
In later sections we make specific reference to the performance of our algorithm when $\rho_{\text{S}} =\rho_{\text{NS}} = 0.75$ as Mathieson \textit{et al.}\cite{Mathieson2016_Performance} set their equivalent parameters to this value in accordance with standards set by clinicians. 

\subsubsection*{Constructing ROC curves}

The important quantities we obtain from our classification procedure are the total number of TPs, FPs, TNs, and FNs in a given recording session which we denote by $N_{\text{TP}}$, $N_{\text{FP}}$, $N_{\text{TN}}$, and $N_{\text{FN}}$, respectively. These are used to define the following quantities which are needed to construct ROC curves:
\\
\textbf{True positive rate (TPR):} $N_{\text{TP}} / (N_{\text{TP}} + N_{\text{FN}}) \in \left[ 0, 1 \right]$, also referred to as the `sensitivity', quantifies the agreement between the expert and the algorithm on detecting the presence of seizure time intervals. 
\\
\textbf{False positive rate (FPR):} $N_{\text{FP}} / (N_{\text{FP}} + N_{\text{TN}}) \in \left[ 0, 1 \right]$, also referred to as `fall-out' or `1 - specificity', quantifies the disagreement between the expert and the algorithm on detecting the presence of non-seizure time intervals. 
\\
\\
Based on the above, the closer the TPR is to $1$ and FPR is to $0$ the better the performance of our algorithm. It is convenient to show this in a two-dimensional plane with FPR values assigned to the x-axis and TPR values to the y-axis. 
In our case, we construct ROC curves by connecting (FPR, TPR) points between $(0,0)$ and $(1,1)$ for increasing values of the FPR where each point corresponds to different choices of the upper threshold $\alpha$ while the other algorithm parameters remain fixed.
Thus, the closer a given point on the ROC curve is to (0,1), the better the algorithm performs. 
We refer to the (0,1) point as the `point of ideal performance'.
\\
\textbf{Why vary $\alpha$ to construct ROC curves: }
We vary the upper threshold $\alpha$ because it is the parameter the algorithm is most sensitive to. Across different GAERS there is a much greater difference in terms of the amplitude of fluctuation than the temporal characteristics of the voltage recordings; seen in part \hyperref[M1]{M1} of the Methods section by comparing the significant differences in voltage amplitudes across GAERS in Fig.~\ref{fig:Morphologies_Rat_T_K_S_} to the marginal differences in the probability densities of seizure and non-seizure time interval durations across GAERS in Fig.~\ref{fig:ResTimes_Rat_T_K_S_}. 
\\
\textbf{Removing misleading points from ROC curves: }
We set TPR$~=0$ if $N_{\text{TP}} + N_{\text{FN}} = 0$ and set FPR$~=0$ if $N_{\text{FP}} + N_{\text{TN}} = 0$ in order to avoid dividing by $0$. 
This typically occurs when  $\alpha$ is not large enough.
See part \hyperref[M4]{M4} of the Methods sections for details on additional steps taken to remove misleading points from ROC curves.
\\
\textbf{Remark:} 
The FPR is related to the `true negative rate' (TNR) via $\text{FPR} = 1- \text{TNR}$, where $\text{TNR} = N_{\text{TN}} / (N_{\text{TN}} + N_{\text{FP}})$.
Temko \textit{et al.} \cite{Temko2011_Performance} and Mathieson \textit{et al.} \cite{Mathieson2016_Performance} construct their ROC curves using the TNR as opposed to the FPR. 
Our preference is to use the FPR, and the difference from using the TNR is a different orientation of the ROC curve.

\subsubsection*{Algorithm performance metrics}

We use the following two metrics to quantify the algorithm performance's in terms of the values of the upper threshold $\alpha$ that are used to construct ROC curves according to the method specified above.
\\
\textbf{Distance from ideal (DFI): }
We quantify the algorithm's performance at individual $\alpha$ values by computing the distance between the corresponding point on the ROC curve and the point of ideal performance, i.e., the (0,1) point. 
We refer to this distance as the `distance from the ideal for a given $\alpha$' and denote it by $\text{DFI}(\alpha)\in[0,\sqrt{2}]$. 
When $\text{DFI}(\alpha) = 0$ this corresponds to ideal performance, meaning all intervals are correctly classified, and $\text{DFI}(\alpha) = \sqrt{2}$ corresponds to the worst possible performance, meaning all intervals are incorrectly classified. 
We say the algorithm achieves excellent performance when $0 < \text{DFI}(\alpha) \leq 0.1$, very good performance when $0.1 < \text{DFI}(\alpha) \leq 0.2$, good performance when $0.2 < \text{DFI}(\alpha) \leq 0.3$, fair performance when $0.3 < \text{DFI}(\alpha) \leq 0.4$, and poor performance when $ \text{DFI}(\alpha) > 0.4$.
We denote the $\alpha$ that minimises $\text{DFI}(\alpha)$ for a given ROC curve as $\alpha^{*}$.
\\
\textbf{Area under the ROC curve (AUC): }
We quantify the algorithm's performance across the range of $\alpha$ values used to construct an ROC curve by computing the `area under the ROC curve'.
We denote this quantity by $\text{AUC} \in [0,1]$. 
When $\text{AUC}=1$, this means the algorithm achieved ideal performance for at least one value of $\alpha$ used to construct the ROC curve.
We say the algorithm achieves excellent performance when $0.9 \leq \text{AUC} < 1$, very good performance when $0.8 \leq \text{AUC} < 0.9$, good performance when $0.7 \leq \text{AUC} < 0.8$, fair performance when $0.6 \leq \text{AUC} < 0.7$, and poor performance when $\text{DFI}(\alpha) < 0.6$. An $\text{AUC} < 0.5$ means the algorithm performed worse than flipping a coin each time a seizure time interval is detected as to whether it will be classified as TP or FP. 
\\
\textbf{Remark: } 
We compute both metrics to provide a more comprehensive assessment of the algorithm's performance.
Further, relying solely on the widely used $\text{AUC}$ metric may result in a misleading evaluation of the algorithm's performance. For instance, when comparing two ROC curves, one may have a larger $\text{AUC}$ but the other may have a smaller $\text{DFI}(\alpha^{*})$, meaning the former shows better performance across a range of $\alpha$ values but the latter shows better performance at a specific $\alpha$ value.
\\
\\
In our experiments, we compare the $\text{DFI}(\alpha^{*})$ obtained from ROC curves that were constructed for the same time series of voltage recordings and range of $\alpha$ values but using different values for the other algorithm parameters. We refer to the set of algorithm parameters corresponding to the smallest $\text{DFI}(\alpha^{*})$ as the `optimal set of parameters' for detecting CTs between NS and S states in a given time series of voltage recordings.


\subsection*{Evaluating the algorithm's performance in a single recording session}

In this subsection we use the $\text{DFI}(\alpha)$ and $\text{AUC}$ metrics to evaluate the algorithm's performance on a single time series of voltage recordings containing multiple seizure and non-seizure time intervals. Through a series of experiments, we highlight the insights these metrics provide and use them to identify the optimal set of algorithm parameters for this time series.

\subsubsection*{Technical details of experiments}

We conduct our experiments using recording session T2M (recordings from rat T during the M$^{th}$ recording session on day 2 of recording). This time series consists of $\approx 7,500\text{s}$ ($\approx 2$ hours) of continuous voltage recordings taken in steps of $0.001\text{s}$, contains 154 seizure intervals according to the expert, i.e., $N_{\text{S}}^{(E)}=154$, and the length of the largest $I_{\text{S}}^{(E)}$ was $\approx 56\text{s}$ and $I_{\text{NS}}^{(E)}$ was $\approx 1200\text{s}$.

We apply our algorithm to the entire T2M time series and classify the resulting $I_{\text{S}}^{(A)}$ and $I_{\text{NS}}^{(A)}$ in each parameter setting specified in Table~\ref{tab:AUC_params}.
More specifically, we consider three different settings of $\tau_{\text{NS}}, ~\tau_{\text{S}}$, and $\tau_{w}$ and refer to these parameter settings as $P_1$, $P_2$, and $P_3$. 
$P_1$ is chosen to reflect the expert's annotation criteria as the length of the shortest $I_{\text{NS}}^{(E)}$ was $\approx 0.7\text{s}$ and shortest $I_{\text{S}}^{(E)}$ was $\approx 1\text{s}$. 
$P_2$ and $P_3$ are chosen to be slightly larger to provide a comparison. 
Furthermore, for $P_1$, $P_2$, and $P_3$, we consider a fixed separation between $\alpha$ and $\beta$, denoted by $\delta_{\alpha \beta} \geq 0$, such that $\beta = \alpha - \delta_{\alpha \beta}$ and $\delta_{\alpha \beta} < \beta \leq \alpha$. We consider values of $\delta_{\alpha \beta} = \left[0.01, 0.06\right]$ and $\alpha \in [\delta_{\alpha \beta}+0.01, 0.1]$. 
We also consider values of $\rho_{\text{S}}, ~\rho_{\text{NS}} \in \left[0, 1 \right]$. 

After using our algorithm to detect CTs in the time series and classifying the resulting $I_{\text{S}}^{(A)}$ and $I_{\text{NS}}^{(A)}$, we construct ROC curves based on TPR and FPR values corresponding to $\alpha \in [\delta_{\alpha \beta}+0.01, 0.1]$ in each of the parameter settings specified in Table~\ref{tab:AUC_params}.

\subsubsection*{Insight from ROC curves - \texorpdfstring{$\rho_{\mathrm{S}}= \rho_{\mathrm{NS}}=0$, $\delta_{\alpha \beta} = 0.01$}{TEXT}}

\begin{figure}[t]
    \centering
    \includegraphics[width=.95\textwidth]{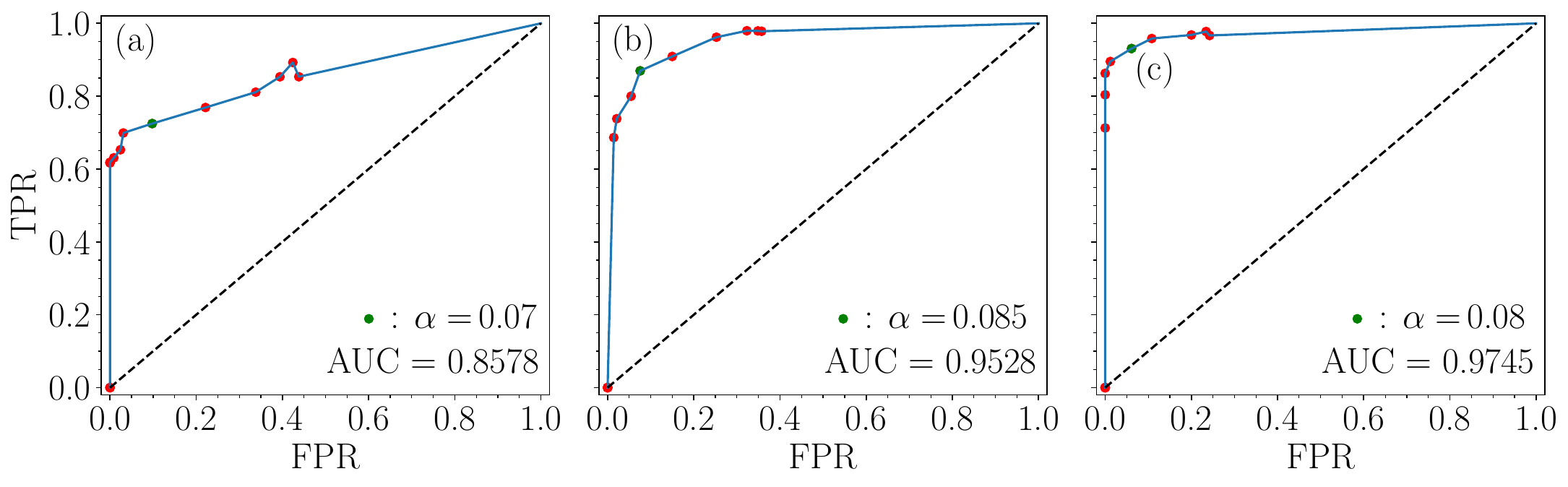}
    \caption{ROC curves obtained when using parameter setting $P_{1}$ in (a), $P_{2}$ in (b), and $P_{3}$ in (c) (see Table~\ref{tab:AUC_params}). In each panel, (FPR, TPR) points plotted in red correspond to different values of $\alpha$, the green point corresponds to $\alpha^{*}$ (the optimal $\alpha$), the lower right corner specifies $\alpha^{*}$ and the AUC.}
    \label{fig:h2_best_rhoS_ROC_example1}
\end{figure}

In Fig.~\ref{fig:h2_best_rhoS_ROC_example1} we show the ROC curves obtained for $P_{1}$ in (a), $P_{2}$ in (b), and $P_{3}$ in (c). In each case $\rho_{\text{S}} = \rho_{\text{NS}}=0$ and $\delta_{\alpha \beta} = 0.01$, meaning $\alpha \in [0.02,0.1]$.
The red data points correspond to (FPR, TPR) values obtained for corresponding values of $\alpha$. 
The green data points correspond to the $\alpha^{*}$ obtained in each case, whose value is specified in the lower right corner. The AUC is also specified in the lower right corner. 

Figure~\ref{fig:h2_best_rhoS_ROC_example1} shows the algorithm performs best in parameter setting $P_{3}$, achieving near ideal performance at several values of $\alpha$. 
More specifically, while the algorithm performs reasonably well in parameter setting $P_{1}$ with $\text{AUC} \approx 0.86$, 
there is a significant improvement in the algorithm's performance for $P_{2}$ with $\text{AUC} \approx 0.95$ and further for $P_{3}$ with $\text{AUC} \approx 0.975$. 
The $\text{DFI}(\alpha^{*})$ for each parameter setting provides similar insight with $\text{DFI}(\alpha^{*}) \approx 0.29$ for $P_{1}$, $\text{DFI}(\alpha^{*}) \approx 0.15$ for $P_{2}$, and $\text{DFI}(\alpha^{*}) \approx 0.09$ for $P_{3}$. 
We found that increasing the values of $\tau_{\text{NS}}$ and $\tau_{\text{S}}$ beyond $P_3$ resulted in little to no improvement in performance but significantly less seizure time intervals were detected by the algorithm in comparison to the expert; see Fig.~\ref{fig:_wider_study_tauNS_tauS_} in [\citeonline{Flynn25_DetSupplement}] for further details.  
\\
\\
Note, the constraint to construct ROC curves for increasing FPR values ensures the AUC can be computed. However, this means that two successive points on a given ROC curve may not necessarily correspond to a successive change in $\alpha$; see Figs.~\ref{fig:P1P2P3_TPR_FPR_alpha_}-\ref{fig:NS_alpha_vs_Ns_P123_deltaAB_0.01_T2M_} in [\citeonline{Flynn25_DetSupplement}] for further details.
\\
\\
In Figs.~\ref{fig:t1_vs_tau1_comparison_vs_alpha_P1_P2_P3_} and \ref{fig:ResTimes_CianAlg} in [\citeonline{Flynn25_DetSupplement}] we present results from studies on the level of agreement achieved between the algorithm and the expert which are not accounted for by ROC curves. Specifically, the agreement between values of $t_{1}$ and $\tau_{1}$.
\\
\\
Further, Temko \textit{et al.} \cite{Temko2011_Performance} and Mathieson \textit{et al.} \cite{Mathieson2016_Performance} use additional metrics derived from their classification procedure to construct curves such as, `seizure detection rate vs. false detections per hour'. We construct similar curves and find they provide complementary insight to ROC curves; see Fig.~\ref{fig:h2_best_rhoS_PV_example1} in [\citeonline{Flynn25_DetSupplement}].


\subsubsection*{AUC for different choices of \texorpdfstring{$\rho_{\mathrm{S}}$}{TEXT} and \texorpdfstring{$\rho_{\mathrm{NS}}$}{TEXT}}

\begin{figure}[t]
    \centering
    \includegraphics[width=.95\textwidth]{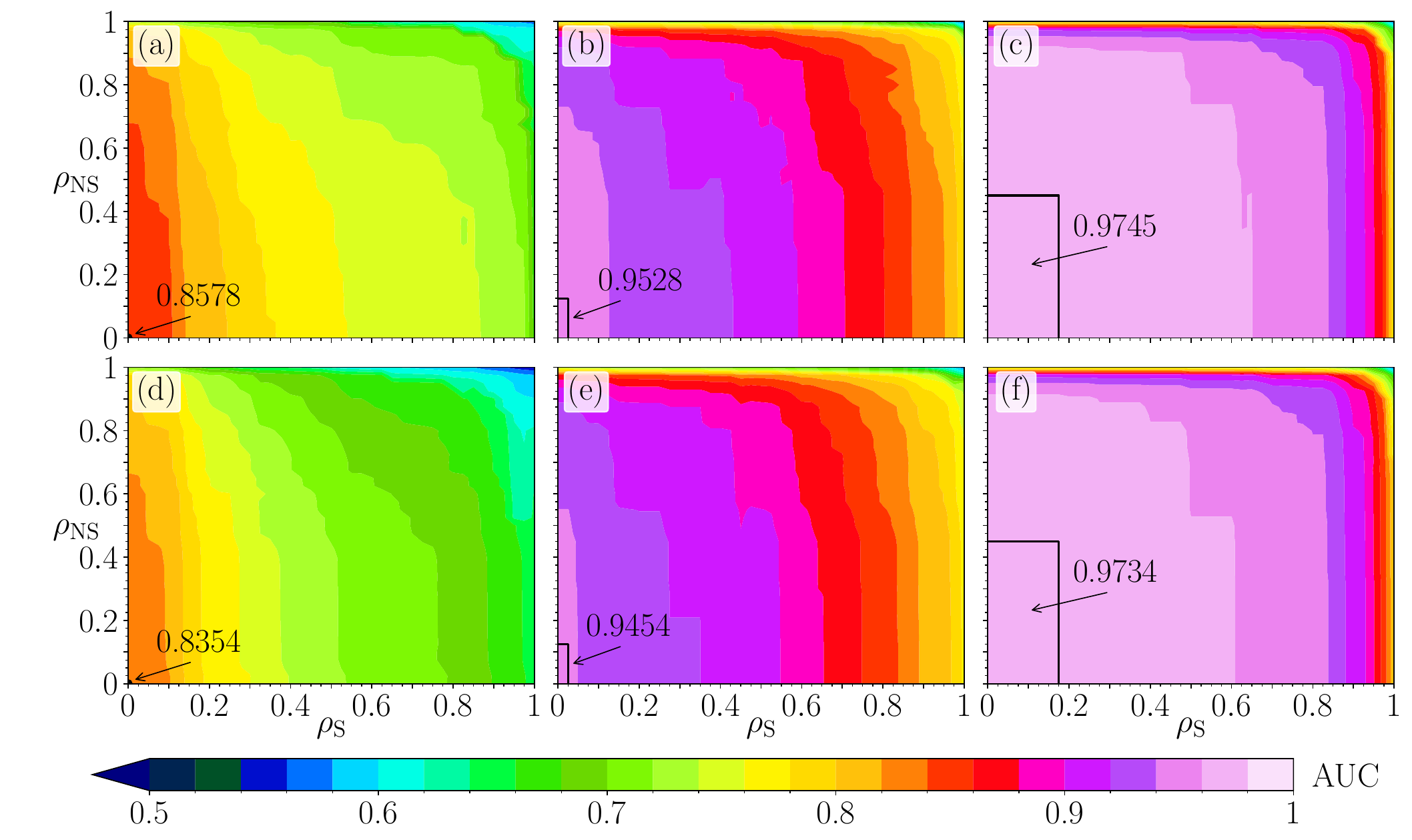}
    \caption{Heatmaps which show the AUC of ROC curves for a given choice of $\rho_{\text{S}}$ and $\rho_{\text{NS}}$. Each column corresponds to parameter settings $P_1$, $P_2$, and $P_3$ as specified in Table~\ref{tab:AUC_params}. Each row corresponds to a different choice of $\delta_{\alpha \beta}$ where $\delta_{\alpha \beta}=0.01$ in (a)-(c) and $\delta_{\alpha \beta}=0.06$ in (d)-(f). The maximum AUC value is highlighted in the lower left corner in (a)-(f).} 
    \label{fig:h2_rhoS_ROC_AUC_}
\end{figure}

In Fig.~\ref{fig:h2_rhoS_ROC_AUC_} we extend the results shown in Fig.~\ref{fig:h2_best_rhoS_ROC_example1} by plotting heatmaps of the AUC values obtained when using different choices of $\rho_{\text{S}}$ and $\rho_{\text{NS}}$ ranging from $0$ to $1$. 
More specifically, the top row of Fig.~\ref{fig:h2_rhoS_ROC_AUC_} shows how the AUC varies for different choices of $\rho_{\text{S}}$ and $\rho_{\text{NS}}$ when using parameter setting $P_1$ in (a), $P_2$ in (b), $P_3$ in (c), where $\delta_{\alpha \beta} = 0.01$ in each case. 
The same is shown in the lower row of Fig.~\ref{fig:h2_rhoS_ROC_AUC_} for $\delta_{\alpha \beta} = 0.06$. 
The black box/point in each panel of Fig.~\ref{fig:h2_rhoS_ROC_AUC_} corresponds to the largest AUC value (up to a difference of $0.001\%$) obtained for a given choice of $\rho_{\text{S}}$ and $\rho_{\text{NS}}$.
There is a box as opposed to a point in some cases as a small number of seizure and non-seizure time intervals are incorrectly classified across a range of $\rho_{\text{S}}$ and $\rho_{\text{NS}}$ values.

Similar to Fig.~\ref{fig:h2_best_rhoS_ROC_example1}, Fig.~\ref{fig:h2_rhoS_ROC_AUC_} shows performance improves as we move the parameter settings away from $P_1$ to $P_2$ and $P_3$. Additionally, Fig.~\ref{fig:h2_rhoS_ROC_AUC_} shows performance worsens as $\delta_{\alpha \beta}$ increases. 
For instance, in Fig.~\ref{fig:h2_rhoS_ROC_AUC_}~(c) ($P_3$ with $\delta_{\alpha \beta} = 0.01$) the AUC is greater than $0.96$ for a wider range of $\rho_{\text{S}}$ and $\rho_{\text{NS}}$ values than in Fig.~\ref{fig:h2_rhoS_ROC_AUC_}~(f) ($P_3$ with $\delta_{\alpha \beta} = 0.06$).
This effect is more pronounced for $P_{1}$ and $P_{2}$.
Furthermore, in terms of the more clinically relevant choice of $\rho_{\text{S}}$ and $\rho_{\text{NS}}$ used by Mathieson \textit{et al.}\cite{Mathieson2016_Performance}, specifically, $\rho_{\text{S}} = \rho_{\text{NS}}=0.75$, the algorithm maintains excellent performance when using $P_{3}$, achieving an $\text{AUC} \approx 0.94$.

We deduce from Figs.~\ref{fig:h2_best_rhoS_ROC_example1} and \ref{fig:h2_rhoS_ROC_AUC_} that the optimal parameter setting of the algorithm for this time series is parameter setting $P_3$ with $\delta_{\alpha \beta} = 0.01$ and $\alpha=0.08$.

\subsection*{Evaluating the algorithm's performance across multiple different recording sessions}

In this subsection we use the AUC and $\text{DFI}(\alpha^{*})$ metrics to evaluate the performance of our algorithm across multiple recording sessions from different GAERS. 
Informed by the results of the previous subsection, we restrict our study to parameter setting $P_3$ with $\delta_{\alpha \beta} = 0.01$. 
Furthermore, we only consider recording sessions where $N_{\text{S}}^{(E)}>60$ in order to base our claims on a reasonably sized sample; see the left column of Table~\ref{tab:ROC_AUC_many_rats_rho_0_075} for details.

\begin{figure}[t]
    \centering
    \includegraphics[width=.95\textwidth]{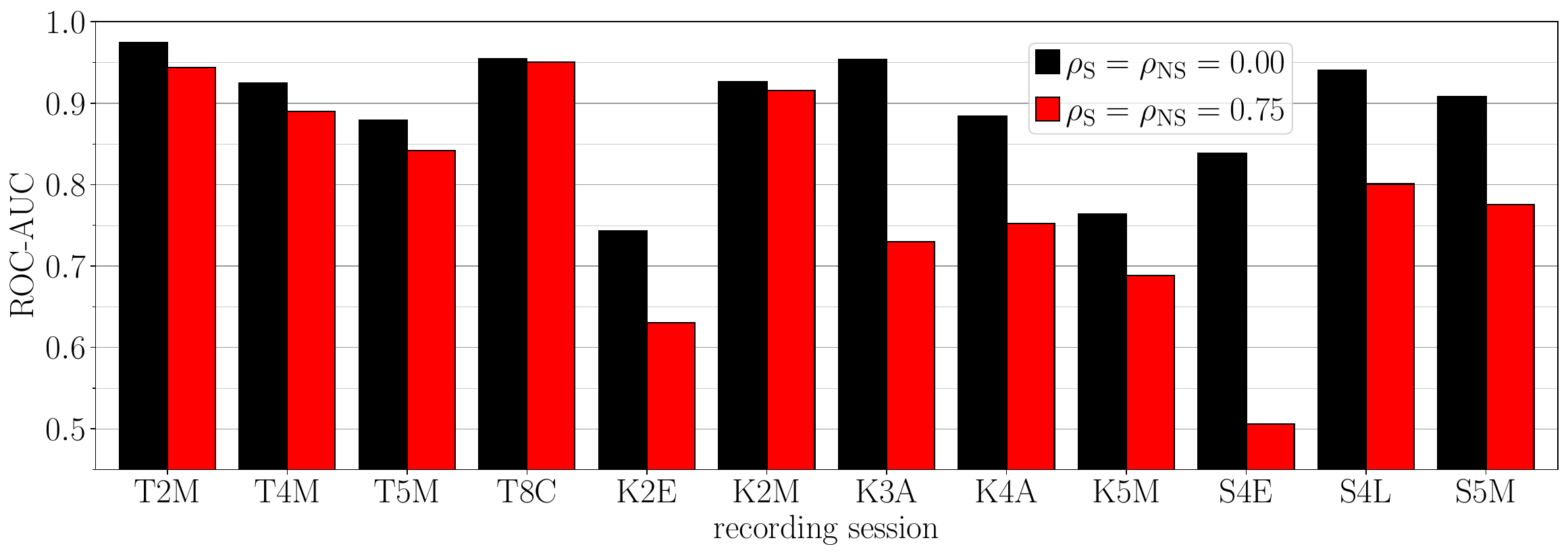}
    \caption{AUC of ROC curves constructed from different recording sessions when (in black) $\rho_{\text{S}} = \rho_{\text{NS}} = 0$ and (in red) $\rho_{\text{S}} = \rho_{\text{NS}} = 0.75$. Remainder of algorithm parameters are chosen as $P_3$ and $\delta_{\alpha \beta}=0.01$ (see Table~\ref{tab:AUC_params}).} 
    \label{fig:_ROC_AUC_multiple_recording_sesseions_rho_0_rho_075_}
\end{figure}

The bar chart in Fig.~\ref{fig:_ROC_AUC_multiple_recording_sesseions_rho_0_rho_075_} shows the resulting AUC for each recording session when choosing $\rho_{\text{S}} = \rho_{\text{NS}} = 0$ (in black) and $\rho_{\text{S}} = \rho_{\text{NS}} = 0.75$ (in red).
In general, Fig.~\ref{fig:_ROC_AUC_multiple_recording_sesseions_rho_0_rho_075_} shows that performance varies depending on which recording session the algorithm is applied to, and in some cases there can be a significant decrease in AUC when $\rho_{\text{S}} = \rho_{\text{NS}} = 0.75$ as opposed to $0$.
More specifically, when $\rho_{\text{S}} = \rho_{\text{NS}} = 0$ the algorithm performs reasonably well, achieving an average $\text{AUC}\approx 0.89$ across the different recording sessions with the largest (best) $\text{AUC}\approx 0.97$ for recording session T2M and smallest (worst) $\text{AUC}\approx 0.75$ for recording session K2E. 
The algorithm still performs reasonably well when $\rho_{\text{S}} = \rho_{\text{NS}} = 0.75$, however, the average AUC decreases to $\approx 0.785$ with the largest (best) $\text{AUC}\approx 0.95$ for recording session T8C and smallest (worst) $\text{AUC}\approx 0.51$ for recording session S4E. 

\begin{figure}[t]
    \centering
    \includegraphics[width=.95\textwidth]{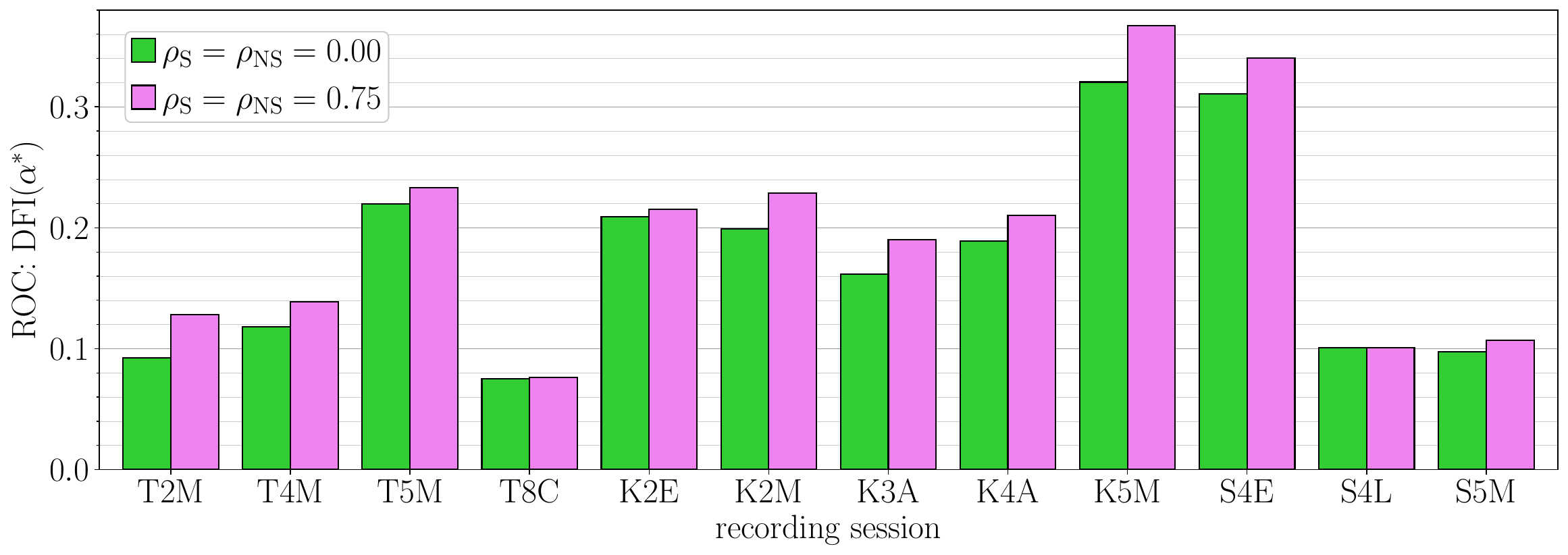}
    \caption{$\text{DFI}(\alpha^{*})$ of ROC curves constructed from different recording sessions when (in green) $\rho_{\text{S}} = \rho_{\text{NS}} = 0$ and (in purple) $\rho_{\text{S}} = \rho_{\text{NS}} = 0.75$. Remainder of algorithm parameters are chosen as $P_3$ and $\delta_{\alpha \beta}=0.01$ (see Table~\ref{tab:AUC_params}).} 
    \label{fig:_DFI_ROC_and_alphastar_multiple_recording_sessions_rho_0_075_}
\end{figure}

The bar chart in Fig.~\ref{fig:_DFI_ROC_and_alphastar_multiple_recording_sessions_rho_0_075_} shows the resulting $\text{DFI}(\alpha^{*})$ for each recording session when choosing $\rho_{\text{S}} = \rho_{\text{NS}} = 0$ (in green) and $\rho_{\text{S}} = \rho_{\text{NS}} = 0.75$ (in purple). Similar to Fig.~\ref{fig:_ROC_AUC_multiple_recording_sesseions_rho_0_rho_075_},  Fig.~\ref{fig:_DFI_ROC_and_alphastar_multiple_recording_sessions_rho_0_075_} shows that the $\text{DFI}(\alpha^{*})$ also varies depending on which recording session the algorithm is applied to.
In contrast to the AUC results, there is generally a relatively smaller change in $\text{DFI}(\alpha^{*})$ when changing $\rho_{\text{S}}$ and $\rho_{\text{NS}}$ from $0$ to $0.75$. 
More specifically, the average $\text{DFI}(\alpha^{*})$ increases from $\approx 0.174$ to $\approx 0.194$, the smallest (best) $\text{DFI}(\alpha^{*})$ increases from $\approx 0.075$ to $\approx 0.076$ for recording session T8C (best in both cases), and the largest (worst) $\text{DFI}(\alpha^{*})$ increases from $\approx 0.32$ to $\approx 0.367$ for recording session K5M (worst in both cases). 

Figures~\ref{fig:_ROC_AUC_multiple_recording_sesseions_rho_0_rho_075_} and \ref{fig:_DFI_ROC_and_alphastar_multiple_recording_sessions_rho_0_075_} show the AUC and $\text{DFI}(\alpha^{*})$ metrics are generally in agreement, i.e., when AUC is relatively large $\text{DFI}(\alpha^{*})$ is relatively small. However, there are some exceptions. For instance, for recording sessions K2E and K2M, the $\text{DFI}(\alpha^{*})$ metric indicates the algorithm achieves very good performance when choosing $\rho_{\text{S}} = \rho_{\text{NS}} = 0$ or $0.75$, but the AUC differs significantly for these recording sessions; the AUC for K2M indicates the algorithm achieves excellent performance, but the AUC for K2E indicates the algorithm achieves only good to fair performance depending on the values of $\rho_{\text{S}}$ and $\rho_{\text{NS}}$.

A more detailed breakdown of the results discussed above is provided in Table~\ref{tab:ROC_AUC_many_rats_rho_0_075}.


\begin{table}[t]
\centering
\def\arraystretch{1.5}
\begin{tabular}{|c|c|c|}
\hline
Rec. sess. ($N_{\text{S}}^{(E)}$) & AUC; $\rho_{\text{S}} = \rho_{\text{NS}}=0.0$ ($\alpha^{*}$: $N_{\text{S}}^{(A)}$, $N_{\text{TP}}$, $N_{\text{TN}}$, DFI) & AUC; $\rho_{\text{S}} = \rho_{\text{NS}}=0.75$ ($\alpha^{*}$: $N_{\text{S}}^{(A)}$, $N_{\text{TP}}$, $N_{\text{TN}}$, DFI) \\
\thickhline
T2M (154) & 0.9745 (0.08: 100, 94, 92, 0.0925) & 0.9440 (0.08: 100, 91, 90, 0.1279) \\
\hline
T4M (135) & 0.9247 (0.07: 101, 99, 87, 0.1182) & 0.8899 (0.07: 101, 96, 86, 0.1386) \\
\hline
T5M (82) & 0.8793 (0.085: 52, 47, 40, 0.2198) & 0.8421 (0.085: 52, 46, 40, 0.2329) \\
\hline
T8C (194) & 0.9545 (0.085: 150, 149, 137, 0.0749) & 0.9508 (0.085: 150, 148, 137, 0.0764) \\
\thickhline
K2E (95) & 0.7435 (0.075: 81, 62, 76, 0.2090) & 0.6303 (0.075: 81, 62, 75, 0.2155) \\
\hline
K2M (112) & 0.9266 (0.08: 118, 106, 96, 0.1992) & 0.9155 (0.08: 118, 105, 92, 0.2287) \\
\hline
K3A (83) & 0.9541 (0.075: 75, 66, 66, 0.1615) & 0.7301 (0.075: 75, 64, 65, 0.1901) \\
\hline
K4A (77) & 0.8841 (0.07: 58, 47, 53, 0.1889) & 0.7521 (0.07: 58, 47, 51, 0.2105) \\
\hline
K5M (244) & 0.7642 (0.085: 124, 94, 97, 0.3205) & 0.6888 (0.08: 125, 86, 98, 0.3673) \\
\thickhline
S4E (127) & 0.8387 (0.07: 69, 62, 44, 0.3110) & 0.5057 (0.06: 98, 71, 77, 0.3402) \\
\hline
S4L (64) & 0.9406 (0.06: 47, 46, 41, 0.1009) & 0.8008 (0.06: 47, 46, 41, 0.1009) \\
\hline
S5M (144) & 0.9078 (0.06: 102, 102, 90, 0.0973) & 0.7753 (0.06: 102, 101, 89, 0.1068) \\
\hline
\end{tabular}
\caption{\label{tab:ROC_AUC_many_rats_rho_0_075} AUC values computed from ROC curves where $\rho_{\text{S}} = \rho_{\text{NS}}=0.0$ (middle column) and $\rho_{\text{S}} = \rho_{\text{NS}}=0.75$ (right column) for specific recording sessions (left column) using parameter setting $P_{3}$ with $\delta_{\alpha \beta} = 0.01$ (see Table~\ref{tab:AUC_params}). 
Information in parenthesis in left column: the no. of seizure time intervals according to the expert. Information in parenthesis in middle and right columns: the optimal $\alpha$ value (denoted by $\alpha^{*}$), and the following corresponding to $\alpha^{*}$: the no. of seizure intervals according to the algorithm, the no. of these time intervals that were classified as TP, and the $\text{DFI}(\alpha^{*})$ (denoted by DFI). 
}
\end{table}

\subsection*{Evaluating the algorithm's performance using a general set of algorithm parameters across multiple recording sessions}

In this subsection we evaluate the performance of our algorithm across multiple recording sessions from different GAERS using a `general' set of algorithm parameters in each case. 
This general set of parameters is based on the average of the $\alpha^{*}$ values listed in Table~\ref{tab:ROC_AUC_many_rats_rho_0_075}. 
Specifically, the average of the $\alpha^{*}$ values when choosing $\rho_{\text{S}} = \rho_{\text{NS}}=0.0$ or $0.75$ is $\approx 0.075$, we denote this average $\alpha$ by $\overline{\alpha}$. 
Thus, we consider parameter setting $P_{3}$ with $\delta_{\alpha \beta}=0.01$ and $\overline{\alpha}$ as our general set of algorithm parameters. 

\begin{figure}[t]
    \centering
    \includegraphics[width=.95\textwidth]{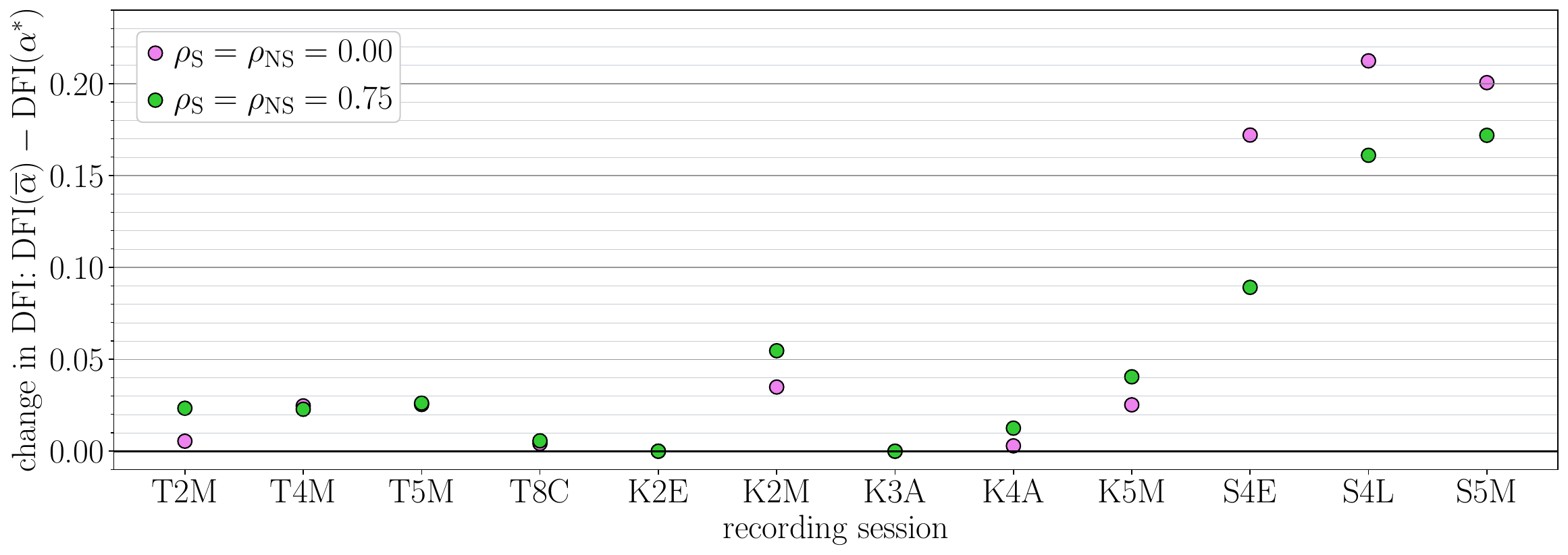}
    \caption{The change in DFI when using an upper threshold of $\overline{\alpha}$ over $\alpha^{*}$ for $\rho_{\text{S}} = \rho_{\text{NS}} = 0$ (in black) and $\rho_{\text{S}} = \rho_{\text{NS}} = 0.75$ (in red). Remainder of algorithm parameters are chosen as $P_3$ and $\delta_{\alpha \beta}=0.01$ (see Table~\ref{tab:AUC_params}).} 
    \label{fig:_changeDFI_alphabar_alphastar_multiple_recording_sessions_rho_0_075_}
\end{figure}

Figure.~\ref{fig:_changeDFI_alphabar_alphastar_multiple_recording_sessions_rho_0_075_} shows the change in performance of our algorithm in terms of the change in the distance from the ideal (DFI) when using the general set of algorithm parameters in comparison to the optimal set of algorithm parameters for each recording session. 
More specifically, we plot $\text{DFI}\left(\overline{\alpha}\right)-\text{DFI}(\alpha^{*})$ for each recording session when choosing $\rho_{\text{S}} = \rho_{\text{NS}}=0.0$ (in green) and $\rho_{\text{S}} = \rho_{\text{NS}}=0.75$ (in purple).
Figure.~\ref{fig:_changeDFI_alphabar_alphastar_multiple_recording_sessions_rho_0_075_} shows that, for many recording sessions, our algorithm maintains much of its performance when using the general set of parameters, i.e., there is little change in the DFI. This is true when $\rho_{\text{S}} = \rho_{\text{NS}} = 0$ or the more clinically relevant setting of $\rho_{\text{S}} = \rho_{\text{NS}} = 0.75$. 
However, for rat S, there is a more noticeable decrease in performance (larger $\text{DFI}(\overline{\alpha})$ than $\text{DFI}(\alpha^{*})$).
This is mostly due to the differences in seizure morphology shown in Fig.~\ref{fig:Morphologies_Rat_T_K_S_}; the amplitudes of  fluctuation in the voltage recordings for rat S are significantly smaller than those for rats T and K. This is further evidenced by the values of $\alpha^{*}$ listed in Table~\ref{tab:ROC_AUC_many_rats_rho_0_075}, i.e., the values of $\alpha^{*}$ for rat S are significantly smaller than those for rats T and K .

\section*{Discussion}

The present paper provides a more detailed analysis on the performance of the seizure detection algorithm introduced in our recent work, Flynn \textit{et al.}~\cite{flynn2025classifying}. 
Specifically, we applied our algorithm to the same set of voltage recordings of seizure activity studied above,
local field potential (LFP) recordings of  seizure activity in \textit{Genetic Absence Epilepsy Rats from Strasbourg} (GAERS), 
and, taking inspiration from Temko \textit{et al.}\cite{Temko2011_Performance} and Mathieson \textit{et al.}\cite{Mathieson2016_Performance}, quantified performance in terms of metrics derived from receiver-operating-characteristic (ROC) curves. 
These metrics show that our algorithm (i) achieves very good performance on average across the voltage recordings we considered when using the optimal set of algorithm parameters in each recording session and (ii) maintains much of its accuracy when using a general set of algorithm parameters applicable across all recording sessions.

While our algorithm only monitors a single feature of the voltage recordings, the voltage at a given time, we expect performance to further improve by monitoring additional features that are relevant to the different types of critical transition (CT) associated with seizure onset, such as the variance and autocorrelation of the time series\cite{flynn2025classifying}. Furthermore, for seizures that involve a bifurcation, it may be possible to detect precursors of seizure onset before a CT by monitoring features associated with `early warning signals of critical slowing down', such as increases in variance and autocorrelation in advance of a CT. However, one must conduct further analysis to determine the predictive power of such features~\cite{lehnertz2024time,ashwin2025_ews_skill}.
To widen our algorithm's applicability, in future work we may consider adaptive threshold techniques to account for seizures whose amplitude either increases or decreases during the course of a seizure~\cite{engel2007epilepsy}.
More generally, our algorithm can be used to detect similar CTs in complex systems beyond the brain, such as CTs found in thermoacoustic systems~\cite{sujith2017thermoacoustic} and aircraft wings~\cite{kurths2025aircraft}.

Overall, our results demonstrate the benefit of studying the epileptic brain through the lens of CTs. We show that seizure detection algorithms can be designed through a mathematics-based approach, without reliance on machine learning tools, and that robust and accurate detection can be achieved despite the presence of artefacts, interictal epileptiform discharges, and variations in seizure morphology.

\section*{Methods}

\subsection*{M1: Obtaining and annotating seizure data from GAERS \label{M1}} 

\begin{figure}[t]
    \centering
    \includegraphics[width=0.85\linewidth]{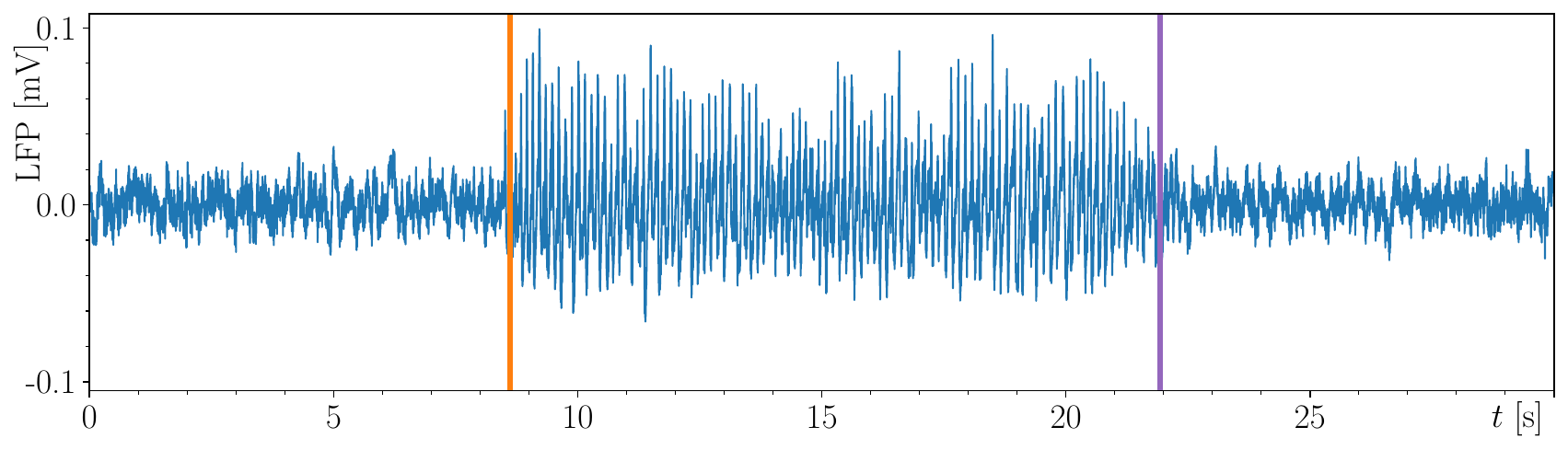}
    \caption{Example of  seizure activity in GAERS. Time series of voltage recordings plotted in terms of local field potential, denoted by LFP vs. time in seconds. 
    Vertical lines correspond to (orange) seizure onset and (purple) seizure offset times according to expert annotations.
    Portion of the voltage recordings shown here is taken from session `S5A'.
    }
    \label{fig:RatData_Properties_}
\end{figure}

\begin{figure}[t]
    \centering
    \includegraphics[width=0.95\linewidth]{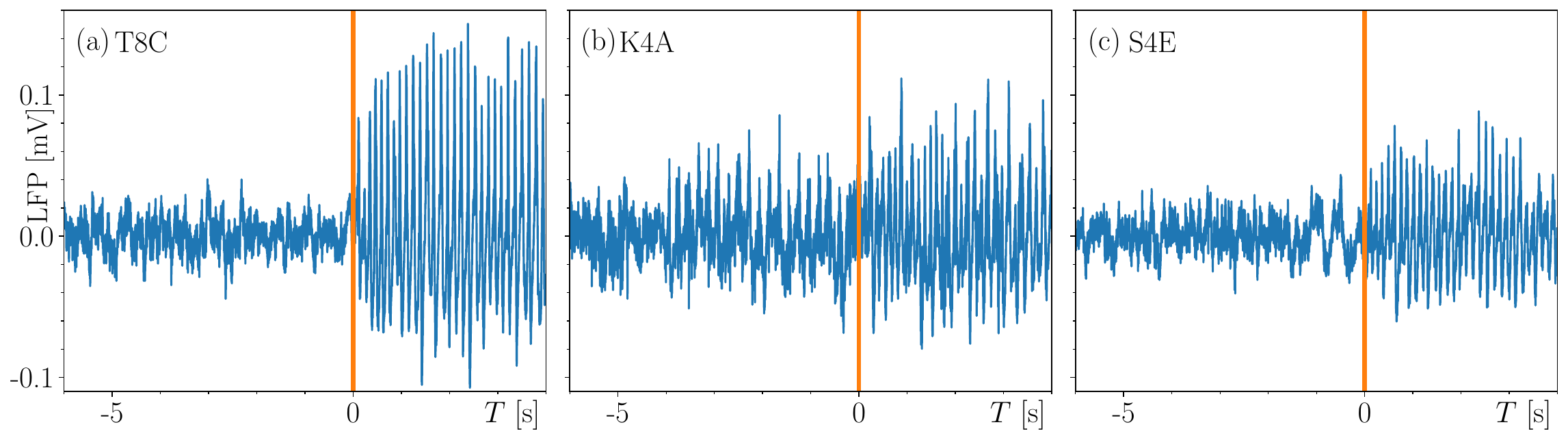}
    \caption{Examples of voltage recordings chosen to highlight some of the differences in seizure morphology between different GAERS (recording session is specified in top left corner). LFP is plotted here vs. $T=t-\tau_{1}$. Vertical orange lines correspond to $T=0$, i.e., the seizure onset time according to the expert.}
    \label{fig:Morphologies_Rat_T_K_S_}
\end{figure}

The data discussed here was obtained by Cian McCafferty, Fran\c{c}ois David, and Vincenzo Crunelli and was presented in \citeonline{mccafferty2018cortical}.
\\
\textbf{Obtaining the data:} The electrophysiological data was acquired from male \textit{Genetic Absence Epilepsy Rats from Strasbourg} (GAERS) aged between 4-7 months when in a state of `relaxed wakefulness' where seizures were experienced more often. Silicon-site electrodes were used to sample voltages (20000/second) from the \textit{ventrobasal thalamus} while the rats were able to move freely and alternate between waking, sleeping, and seizing states. This data was processed with a Plexon HST/32V-G20 VLSI-based preamplifier and associated digitization system and subsequently down-sampled to 1000 samples/second. 
\\
\textbf{Labelling the recording session:}
Each recording session was labelled using the following convention described through example, `S1K': voltage recordings from rat S on day 1 of recording during the K$^{th}$ recording session on that day.
\\
\textbf{Annotating the data:}
We denote seizure onset times by $\tau_{1}$ and offset times by $\tau_{2}$. The process used by the expert to obtain these is specified below.
Spike-wave discharges (SWDs) are the electrical hallmark of absence seizures, they define electrical seizure onset and offset times. These were identified using Cambridge Electronic Design's `Spike2' software and the following procedure. In all cases, EEG at $1000\text{Hz}$ was used for this step. EEG was acquired at $1000\text{Hz}$ for fMRI (see McCafferty \textit{et al.}\cite{mccafferty2023decreased}) and behaviour was down-sampled by averaging for neuronal activity. 
Briefly, smoothening (voltage at time $t$ is set to the mean of voltages from time $t-10\text{ms}$ to time $t+10\text{ms}$) and DC removal (voltage at time $t$ is set to the original voltage minus the mean of voltages from time $t-0.1\text{s}$ to time $t+0.11\text{s}$) functions were used to reversibly visually clean the frontoparietal differential EEG. 
Then, a negative amplitude threshold (mean voltage minus $5$-$7$~standard deviations of baseline non-SWD EEG) was used to detect putative spike-wave crossing points, defined as whenever the signal crossed this amplitude threshold. 
The crossing points were then grouped into events based on the intervals between them (maximum time between initial two crossings $0.2\text{s}$, maximum time between any two crossings within an event $0.35\text{s}$, minimum of 5 crossings per event) and the defined properties of SWDs (minimum duration $0.5\text{s}$, minimum inter-SWD interval $0.5\text{s}$ merging any SWDs with shorter intervals), and subsequently, using a frequency threshold, these events were classified as SWDs (if $>75\%$ of intercrossing intervals were within a 5–12\text{Hz} range) or other (e.g., noise, sleep).
Labelled SWDs were then visually inspected for accuracy of SWD detection, as well as seizure onset and offset times. 
Periods of sleep were identified based on sharp increases in the 1–4\text{Hz} frequency band and were excluded from analysis.
Periods of non-REM sleep were rare in the recordings and not informative due to their relatively short duration.
\\
\textbf{Examples and statistics of seizure activity: }
See Fig.~\ref{fig:RatData_Properties_} for an example of expert annotations of voltage recordings.
See Fig.~\ref{fig:Morphologies_Rat_T_K_S_} for representative examples of seizure activity from different GAERS, illustrating differences in seizure morphology. 
Specifically, the average amplitude of fluctuation in the S state is much larger in Fig.~\ref{fig:Morphologies_Rat_T_K_S_}~(a) than in Figs.~\ref{fig:Morphologies_Rat_T_K_S_}~(b) and (c), a difference that is consistently observed across recording sessions from each of these GAERS. 
While Fig.~\ref{fig:Morphologies_Rat_T_K_S_} shows differences in how the S state emerges, this is not unique to individual GAERS.
We show in Flynn \textit{et al.}~\cite{flynn2025classifying} that the S state can emerge through different types of CT, with similar frequencies of each transition type observed across GAERS.


\begin{figure}[t]
    \centering
    \includegraphics[width=0.999\linewidth]{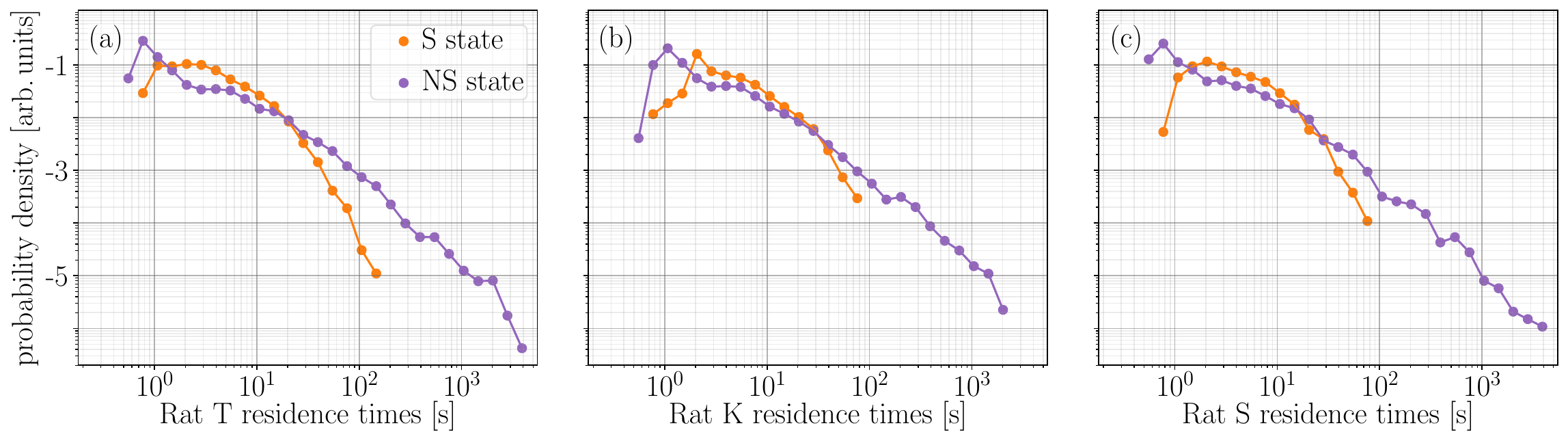}
    \caption{Probability density of residence times in the S and NS states computed from expert annotations of seizure onset and offset times for all voltage recordings of rat T in (a), rat K in (b), and rat S in (c). The information presented here is based on 621, 1593, and 1136 transitions from each respective rat.}
    \label{fig:ResTimes_Rat_T_K_S_}
\end{figure}

See Fig.~\ref{fig:ResTimes_Rat_T_K_S_} for a comparison of the probability density of residence times in the S and NS states for all seizure onset and offset times annotated by the expert in the voltage recordings of rats T, K, and S. Note, there are similar distributions for each rat. The average residence time in the S and NS state for rat T was $\approx11~\text{s}$ and $\approx86~\text{s}$, for rat K was $\approx13~\text{s}$ and $\approx76~\text{s}$, and for rat S was $\approx10~\text{s}$ and $\approx74~\text{s}$.

In summary, when comparing the insights on seizure activity in Fig.~\ref{fig:Morphologies_Rat_T_K_S_} to Fig.~\ref{fig:ResTimes_Rat_T_K_S_}, we observe much greater variability in terms of the amplitude of fluctuation than the temporal characteristics of seizure activity across different GAERS. 

\begin{figure}[t]
    \centering
    \includegraphics[width=0.7\linewidth]{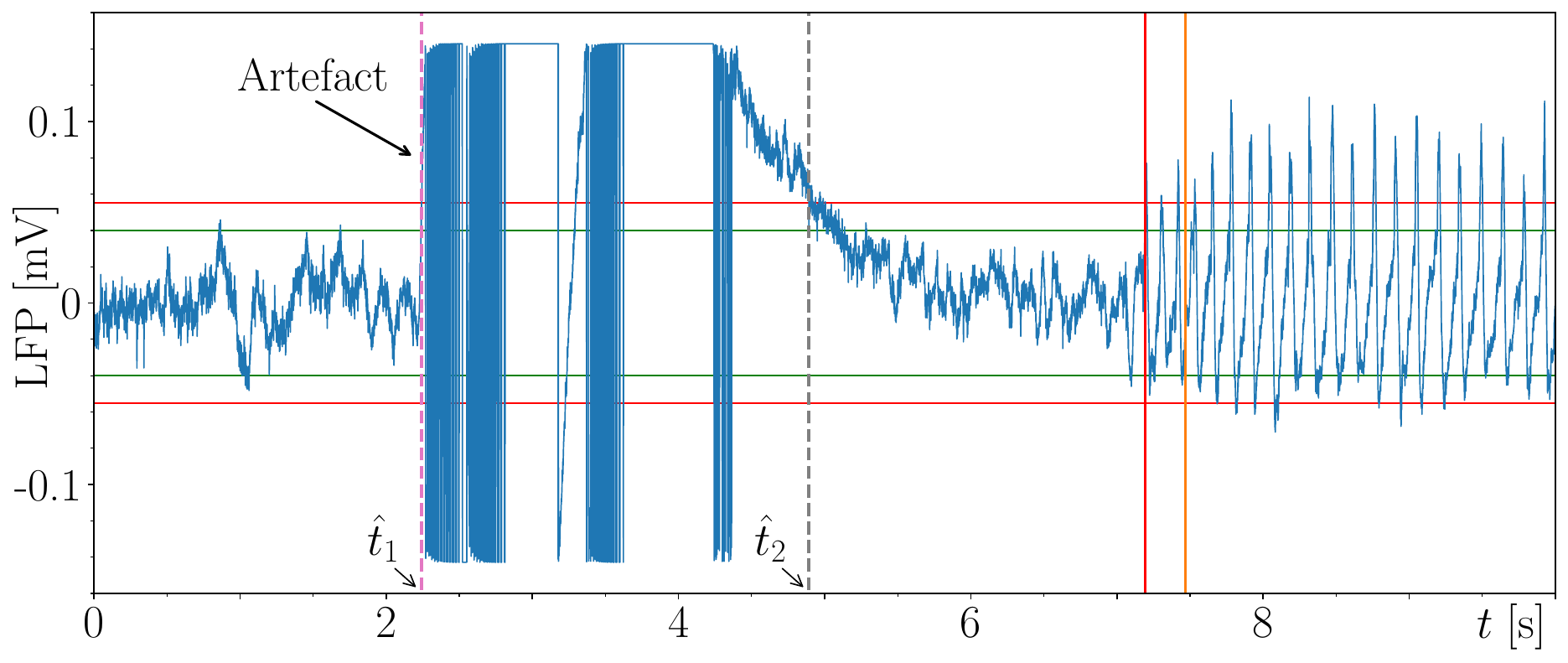}
    \caption{
    Example of artefact activity in voltage recordings. Vertical dashed lines indicate (in pink) the beginning and (in grey) end of artefact activity. Vertical solid lines indicate seizure onset time according to (in orange) the expert and (in red) the algorithm when applied to a portion of the voltage recordings from session `T8C'. 
    Algorithm parameters chosen as $\alpha=0.055$, $\beta=0.04$, $\tau_{\text{NS}}=5$, $\tau_{\text{S}}=2$, $\tau_{w}=1$, and $\Delta = 0.001$. Horizontal lines indicate the thresholds of (in red) $\alpha$ and (in green) $\beta$.}
    \label{fig:Compare_detection_RatData_artefact}
\end{figure}

\noindent
\textbf{Artefacts and how they are accounted for:}
In the voltage recordings we analysed, artefacts are high-amplitude and high-frequency oscillatory-like events; see Fig.~\ref{fig:Compare_detection_RatData_artefact} for an example. 
Artefacts appear for different reasons, the most common being (i) movement of electrodes/wiring and (ii) movement of muscles close to the electrodes (generally for chewing). 
Some artefacts are generated by signal overload and noise from electrical devices around the recording setup, however, these are less common as the recording process was generally appropriately amplified and shielded. 
Artefacts are accounted for through the annotation method described above since these events do not align with the frequency profile of SWDs. 
Additionally, the visual inspection described above also includes the manual rejection of artefacts which were not excised from annotation. 
Entire recording sessions were excluded if more than $5\%$ of the session consisted of artefact activity.

\noindent
\textbf{Remark (comparison with human data):}
GAERS can express multiple absence seizures per minute\cite{powell2014seizure}, far exceeding the frequency of any absence seizure syndrome in humans. 
For instance, Gregor{\v{c}}i{\v{c}} \textit{et al.}\cite{gregorvcivc2022difficult} found that, for children with treatment-resistant childhood absence epilepsy, the median number of seizures per day was three.
For a review of the GAERS model, with attention to its similarities and differences to human absence seizures, see Depaulis \textit{et al.}\cite{depaulis2016genetic}.


\subsection*{M2: CT detection algorithm applied to voltage recordings of  seizure activity - technical details \label{M2}}

\begin{figure}[t]
    \centering
    \includegraphics[width=0.7\linewidth]{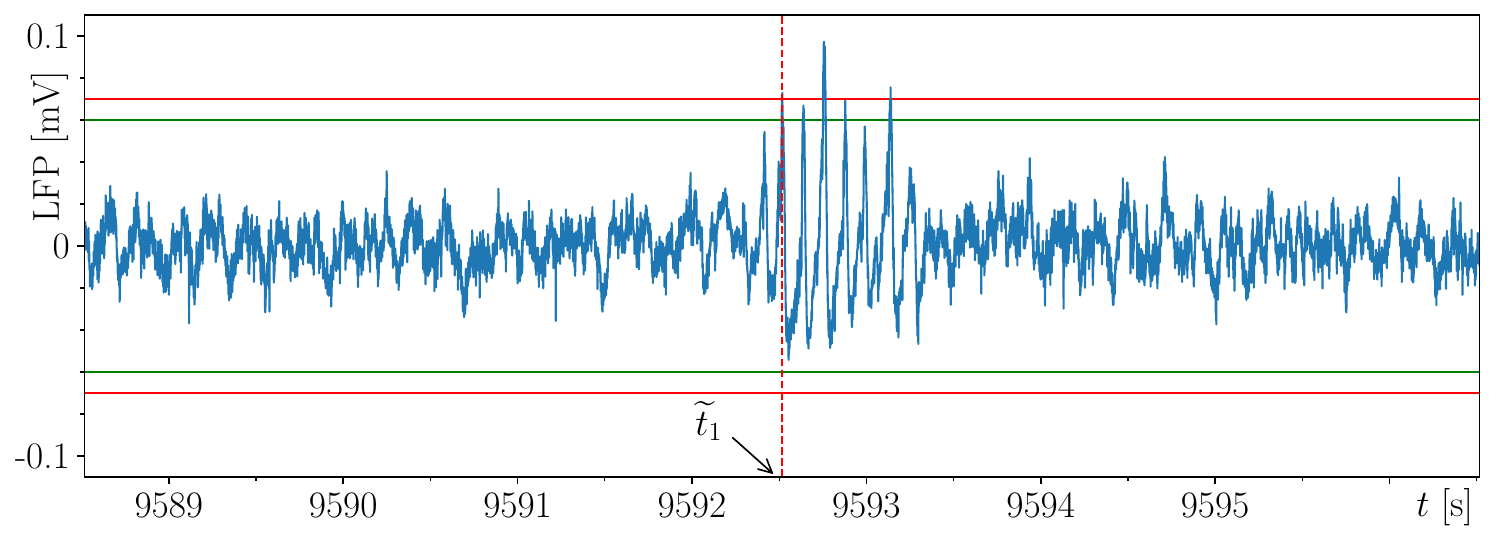}
    \caption{Example of an almost-occurring CT (interictal epileptiform discharges). Vertical dashed line indicates (in red) when a CT from the NS to S state almost occurs at time $t=\tilde{t}_{1}$. 
    Algorithm parameters chosen as $\alpha=0.07$, $\beta=0.06$, $\tau_{\text{NS}}=3$, $\tau_{\text{S}}=2$, $\tau_{w}=1$, and $\Delta = 0.001$. Horizontal lines indicate the thresholds of (in red) $\alpha$ and (in green) $\beta$.}
    \label{fig:AlmostCT_example_}
\end{figure}

For convenience, we describe the algorithm in terms of a variable $x(t)$, which in our case is the LFP at a given time $t$, i.e., $\text{LFP}(t)$. 
We use $\delta > 0$ to denote the time interval between two consecutive data points in a given time series, in the case of voltage recordings in GAERS, we have $\delta=0.001\text{s}$. 
We introduce the following six parameters that our algorithm uses: 
the upper voltage threshold $\alpha > 0$,
the lower voltage threshold $0 <\beta < \alpha$,
the size of the moving window $\tau_w > 0$,
the time step size of the moving window $\delta \leq \Delta \leq \tau_w$,
the minimum time duration of larger-amplitude fluctuations, $\tau_{\text{S}} > \tau_{w}$, expressed as $\tau_{\text{S}} = n_{\text{S}}\Delta + \tau_w$, where $n_\text{S}$ is an integer, 
and the minimum time duration of smaller-amplitude fluctuations, $\tau_{\text{NS}} \geq \tau_{\text{S}}$, expressed similarly as $\tau_{\text{NS}} = n_{\text{NS}}\Delta + \tau_w$, where $n_{\text{NS}}$ is an integer. In all our experiments we set $\Delta = \delta = 0.001$.
\\
{\bf The starting point:} 
We start in the NS state, where $|x(t)| < \beta$ for the time duration of at least $\tau_{\text{NS}}$.
\\
{\bf The moving window:} 
When the brain is in the NS state and $|x(t)|$ exceeds $\alpha$ at time $t = t_j$, the moving window is activated and $|x(t)|$ is examined within consecutive windows of duration $\tau_w$ that are shifted in time by $\Delta$, starting with $[t_j,\, t_j + \tau_w]$, then $[t_j + \Delta,\, t_j  + \tau_w + \Delta]$, $[t_j + 2 \, \Delta,\, t_j  + \tau_w + 2 \, \Delta]$, and so on. 
The moving window is deactivated in two cases: (i) the brain is in the NS state, the window is activated, but the algorithm does not detect a CT to the S state, and (ii) the brain is in the S state and the algorithm detects a CT to the NS state.  
\\
{\bf Critical transitions:}
The algorithm detects a {\em CT from the NS to S state} at time $t=t_1$ if:
\begin{itemize}[labelindent=10pt,itemindent=1em,leftmargin=!]
    \item[(a1)\label{a1}]
    The brain is in the NS state just before $t_1$.
    \item[(a2)\label{a2}]
    $|x(t)|$ exceeds $\alpha$ at time $t=t_1$,  i.e., $|x(t_1)|=\alpha$ and $|x(t_1 + \delta)| > \alpha$.
     \item[(a3)\label{a3}]
     Each of the $n_{\text{S}}$ consecutive positions of the moving window contains an $|x(t)| > \beta$.
\end{itemize}
The algorithm detects a {\em CT from the S to NS state} at time  $t=t_2$ if:
\begin{itemize}[labelindent=10pt,itemindent=1em,leftmargin=!]
    \item[(b1)\label{b1}]
    The brain is in the S state just before $t_2$.
    \item[(b2)\label{b2}]
    $|x(t)|$ falls below $\beta$ at time $t=t_2$, i.e., $|x(t_2)|=\beta$ and $|x(t_2 + \delta)| < \beta$.
     \item[(b3)\label{b3}]
     Each of the $n_{\text{NS}}$ consecutive positions of the moving window contains no $|x(t)| \geq \alpha$.
\end{itemize}
In other words, the algorithm detects a CT from the NS to S state if $|x(t)|$ exceeds the upper threshold $\alpha$ and then continues to exceed the lower threshold $\beta$ frequently enough for a period of at least $\tau_{\text{S}}$. Similarly, the algorithm detects a CT from the S to NS state if $|x(t)|$ falls below the lower threshold $\beta$ and then does not exceed the upper threshold $\alpha$ for a period of at least $\tau_{\text{NS}}$. 
\\
\textbf{Almost-occurring critical transitions: }
The algorithm detects an {\em almost-occurring CT from the NS to S state} at time  $t=\tilde{t}_{1}$ if \hyperref[a1]{(a1)} and \hyperref[a2]{(a2)} are satisfied but \hyperref[a3]{(a3)} is not. 
Similarly, the algorithm detect an {\em almost-occurring CT from the S to NS state} at time  $t=\tilde{t}_{2}$ if \hyperref[b1]{(b1)} and \hyperref[b2]{(b2)} are satisfied but \hyperref[b3]{(b3)} is not.
See Fig.~\ref{fig:AlmostCT_example_} for an example of an almost-occurring CT from the NS to S state. Note, events like these are associated with interictal epileptiform discharges, short time intervals of seizure-like activity that do not meet the criteria to be considered as seizure activity.

\subsubsection*{Accounting for artefacts}

Voltage recordings of  seizure activity are notoriously susceptible to artefacts. The data we analyse in this paper is no exception, see Fig.~\ref{fig:Compare_detection_RatData_artefact} for a typical example of artefact activity and part \hyperref[M1]{M1} of the Methods section for reasons why artefacts appear.
Artefacts can be mistaken for seizure activity and seizure detection algorithms need to be designed accordingly. 
We now describe how we alter the above algorithm to account for artefacts.

From Fig.~\ref{fig:Compare_detection_RatData_artefact} we observe that for $t \in [2,5]$, the LFP jumps between $-0.15$ and $0.15$ mV in short time intervals, often within consecutive measurements ($0.001\text{s}$). 
Based on this observation, we alter condition \hyperref[a2]{(a2)} of our algorithm to account for artefacts as follows: 
\begin{itemize}[labelindent=20pt,itemindent=1em,leftmargin=!]
    \item[(a2-1)\label{a2-1}] $|x(t)|$ exceeds $\alpha$ at time $t=t_{1}'$ \textbf{and} $|x(t+\delta)-x(t)| < \xi$ for all $t \in \left[t_{1}', t_{1}' + t_{w} \right]$, or,
    \item[(a2-2)\label{a2-2}] $|x(t)|$ exceeds $\alpha$ at time $t=t_{1}'$ \textbf{and} $|x(t+\delta)-x(t)| \geq \xi$ for any $t \in \left[t_{1}', t_{1}' + t_{w} \right]$.
\end{itemize}
If \hyperref[a2-1]{(a2-1)} is true then the algorithm continues as before to evaluate whether a CT from the NS to S state is detected at $t_{1}'$. 
On the other hand, if \hyperref[a2-2]{(a2-2)} is true then we say artefact activity begins at time $\hat{t}_{1}=t_{1}'$ and we continue to monitor $|x(t)|$ in moving windows. 
We say artefact activity ends at time $t=\hat{t}_{2}$ if $|x(\hat{t}_{2})| < \alpha$ \textbf{and} $|x(t+\delta)-x(t)| < \xi$ for all $t \in \left[\hat{t}_{2}, \hat{t}_{2} + t_{w} \right]$. 
We find most artefacts are accounted for by choosing $\xi = 0.2$. 

Fig.~\ref{fig:Compare_detection_RatData_artefact} shows the result of applying our algorithm with \hyperref[a2-1]{(a2-1)} and \hyperref[a2-2]{(a2-2)} and parameters chosen as $\alpha=0.055$, $\beta=0.04$, $\tau_{\text{NS}}=5$, $\tau_{\text{S}}=2$, $\tau_{w}=1$, and $\Delta = 0.001$. 
The example in Fig.~\ref{fig:Compare_detection_RatData_artefact} is chosen to show that for the current choice of $\tau_{\text{S}}$ and $\tau_{\text{NS}}$, the time when artefact activity begins would have been considered as the time when a CT from the NS to S state occurs unless artefacts are accounted for. 
Furthermore, Fig.~\ref{fig:Compare_detection_RatData_artefact} shows that by excluding these artefacts there is still a strong agreement between the expert and the algorithm on the seizure onset times. 


\subsection*{M3: Seizure and non-seizure time interval classification procedure \label{M3}}

\subsubsection*{Seizure and non-seizure time intervals}

We introduce the following terminology to define seizure and non-seizure time intervals in time series of voltage recordings based on expert annotations of seizure onset and offset times and the times that our algorithm detects CTs.
\\
\textbf{Expert: }
We denote the seizure onset and offset times according to the expert as $\tau_{1}^{(i)}$ and $\tau_{2}^{(i)}$ where $i=1, 2, \ldots, N_{\text{S}}^{(E)}$ and $N_{\text{S}}^{(E)}$ is the number of seizure time intervals in the time series according to the expert. 
We define the $i^{th}$ seizure time interval according to the expert as $\left(I_{\text{S}}^{(E)}\right)_{i} = \left\{ t \in \mathbb{R} \,:\, \tau_{1}^{(i)} \leq t < \tau_{2}^{(i)} \right\}$. 
We define the subsequent non-seizure time interval as $\left(I_{\text{NS}}^{(E)}\right)_{i} = \left\{ t \in \mathbb{R} \,:\, \tau_{2}^{(i)} \leq t < \tau_{1}^{(i+1)} \right\}$. 
Thus, for a given $N_{\text{S}}^{(E)}$, there are $N_{\text{NS}}^{(E)}= N_{\text{S}}^{(E)}-1$ non-seizure time intervals. 
\\
\textbf{Algorithm: } 
We denote the times that our algorithm detects CTs from the NS to S state as $t_{1}^{(j)}$ and CTs from the S to NS state as $t_{2}^{(j)}$ where $j=1, 2, \ldots, N_{\text{S}}^{(A)}$ and $N_{\text{S}}^{(A)}$ is the number of seizure time intervals in the time series according to the algorithm. 
We use the same convention to define the $j^{th}$ seizure time interval according to the algorithm as $\left(I_{\text{S}}^{(A)}\right)_{j} = \left\{ t \in \mathbb{R} \,:\, t_{1}^{(j)} \leq t < t_{2}^{(j)} \right\}$. 
We define the subsequent non-seizure time interval as $\left(I_{\text{NS}}^{(A)}\right)_{j} = \left\{ t \in \mathbb{R} \,:\, t_{2}^{(j)} \leq t < t_{1}^{(j+1)} \right\}$. 
Similarly, for a given $N_{\text{S}}^{(A)}$, there are $N_{\text{NS}}^{(A)}= N_{\text{S}}^{(A)}-1$ non-seizure time intervals. 
\\
\\
The length of the above time intervals is calculated as follows, $\left|\left(I_{\text{S}}^{(A)}\right)_{j}\right| = t_{2}^{(j)} - t_{1}^{(j)}$, $\left|\left(I_{\text{NS}}^{(A)}\right)_{j}\right| = t_{1}^{(j+1)} - t_{2}^{(j)}$, and similarly for $\left|\left(I_{\text{S}}^{(E)}\right)_{i}\right|$ and $\left|\left(I_{\text{NS}}^{(E)}\right)_{i}\right|$ using the corresponding $\tau_{1}$ and $\tau_{2}$ terms, where $|I|$ denotes the length of the interval $I$. 
We refer to a given $\left|\left(I_{\text{S}}^{(A)}\right)_{j}\right|$ as a residence time in the S state according to the algorithm, $\left|\left(I_{\text{NS}}^{(A)}\right)_{j}\right|$ as a residence time in the NS state, and similarly for $\left|\left(I_{\text{S}}^{(E)}\right)_{i}\right|$ and $\left|\left(I_{\text{NS}}^{(E)}\right)_{i}\right|$ in terms of the expert's annotations.

\subsubsection*{Classification of \texorpdfstring{$A_{\mathrm{S}}$}{TEXT} and \texorpdfstring{$A_{\mathrm{NS}}$}{TEXT}}

The metrics used to evaluate the performance of our algorithm are based on the following time interval classification procedure: 
\\
\textbf{True positive: } we classify a given $\left(I_{\text{S}}^{(A)}\right)_{j}$ as {\em true positive} (TP) if any of the following three conditions are true:
\begin{itemize}\setlength{\itemindent}{.1in}
    \item[C1: \label{C1}] $\left(I_{\text{S}}^{(A)}\right)_{j}$ is contained in some  $\left(I_{\text{S}}^{(E)}\right)_{i}$, that is $\left(I_{\text{S}}^{(A)}\right)_{j}\subseteq \left(I_{\text{S}}^{(E)}\right)_{i}$ for some $i$. 
    \item[C2: \label{C2}] $\left(I_{\text{S}}^{(A)}\right)_{j}$ contains at least one $\left(I_{\text{S}}^{(E)}\right)_{i}$, that is $\left(I_{\text{S}}^{(A)}\right)_{j} \supset \left(I_{\text{S}}^{(E)}\right)_{i}$ for some $i$.
    \item[C3: \label{C3}] {Neither \hyperref[C1]{C1} or \hyperref[C2]{C2} are true but there is sufficient overlap between $\left(I_{\text{S}}^{(A)}\right)_{j}$ and some $\left(I_{\text{S}}^{(E)}\right)_{i}$,\\ 
    \phantom{x}\hspace{.08cm}that is $| \left(I_{\text{S}}^{(A)}\right)_{j} \bigcap \left(I_{\text{S}}^{(E)}\right)_{i} | / |\left(I_{\text{S}}^{(E)}\right)_{i}| > \rho_{\text{S}}$ for some  $i$ and a given $0 \leq \rho_{\text{S}} \leq 1$.} 
\end{itemize}
\textbf{False positive: } if \hyperref[C1]{C1}-\hyperref[C3]{C3} are not true then $\left(I_{\text{S}}^{(A)}\right)_{j}$ is classified as {\em false positive} (FP). 
\\
\textbf{True negative: } we classify a given $\left(I_{\text{NS}}^{(A)}\right)_{j}$ as {\em true negative} (TN) if any of the following three conditions are true: 
\begin{itemize}\setlength{\itemindent}{.1in}
    \item[D1: \label{D1}] $\left(I_{\text{NS}}^{(A)}\right)_{j}$ is contained in some  $\left(I_{\text{NS}}^{(E)}\right)_{i}$, that is $\left(I_{\text{NS}}^{(A)}\right)_{j}\subseteq \left(I_{\text{NS}}^{(E)}\right)_{i}$ for some $i$.
    \item[D2: \label{D2}] $\left(I_{\text{NS}}^{(A)}\right)_{j}$ contains one $\left(I_{\text{NS}}^{(E)}\right)_{i}$ and no $\left(I_{\text{S}}^{(E)}\right)_{k}$, that is $\left(I_{\text{NS}}^{(A)}\right)_{j}\supset \left(I_{\text{NS}}^{(E)}\right)_{i}$ for some $i$. 
    \item[D3: \label{D3}] Neither \hyperref[D1]{D1} or \hyperref[D2]{D2} are true but there is an $\left(I_{\text{NS}}^{(E)}\right)_{i}$ with sufficient overlap with $\left(I_{\text{NS}}^{(A)}\right)_{j}$,\\
    \phantom{x}\hspace{.08cm}that is $| \left(I_{\text{NS}}^{(A)}\right)_{j} \bigcap \left(I_{\text{NS}}^{(E)}\right)_{i} | / |\left(I_{\text{NS}}^{(A)}\right)_{j}| > \rho_{\text{NS}}$ for some  $i$ and a given $0 \leq \rho_{\text{NS}} \leq 1$ \textbf{and} \\
    \phantom{x}\hspace{.08cm}$\left(I_{\text{NS}}^{(A)}\right)_{j}$ does not contain any  $\left(I_{\text{S}}^{(E)}\right)_{k}$, i.e.,  $\left(I_{\text{NS}}^{(A)}\right)_{j} \nsupset  \left(I_{\text{S}}^{(E)}\right)_{i}$ for all $k$. 
\end{itemize}
\textbf{False negative: } if \hyperref[D1]{D1}-\hyperref[D3]{D3} are not true then $\left(I_{\text{NS}}^{(A)}\right)_{j}$ is classified as {\em false negative} (FN). 
\\
\\
We denote the total number of TPs with $N_{\text{TP}}$, FPs with $N_{\text{FP}}$, TNs with $N_{\text{TN}}$, and FNs with $N_{\text{FN}}$. Table~\ref{tab:Case_i_a} specifies how we check \hyperref[C1]{C1}-\hyperref[C3]{C3} and \hyperref[D1]{D1}-\hyperref[D3]{D3} in terms of the corresponding $t_{1}$, $t_{2}$, $\tau_{1}$, and $\tau_{2}$ values. 
\\
\\
Note, in \hyperref[C3]{C3} and \hyperref[D3]{D3} we introduced the parameters, $\rho_{\text{S}}$ and $\rho_{\text{NS}}$, which we use as thresholds to quantify the minimum amount of overlap between time intervals that is needed for a given $I_{\text{S}}^{(A)}$ or $I_{\text{NS}}^{(A)}$ to be classified as TP or TN. These parameters can be chosen in accordance with standards set by clinicians, for instance, in Mathieson \textit{et al.}\cite{Mathieson2016_Performance} the authors set their parameters to $0.75$.  
\\
\\
\textbf{Remark:} based on \hyperref[C2]{C2}, we allow for a given $\left(I_{\text{S}}^{(A)}\right)_{j}$ to be classified as TP even if it contains more than one $\left(I_{\text{S}}^{(E)}\right)_{i}$ and an $\left(I_{\text{NS}}^{(E)}\right)_{k}$. However, according to certain performance metrics, this can result in misleadingly high performance for reasons that we outline in part \hyperref[M4]{M4} of the Methods Section. 
In contrast to \hyperref[C2]{C2}, with \hyperref[D2]{D2} we only allow for a given $\left(I_{\text{NS}}^{(A)}\right)_{j}$ to be classified as TN if it contains no more than one $\left(I_{\text{NS}}^{(E)}\right)_{i}$.

\begin{table}[t]
\centering
\def\arraystretch{1.2}
\begin{tabular}{|c||c|}
\hline
\multirow{2}{4em}{\centering \hyperref[C1]{C1}} & \multirow{2}{24em}{\centering $\tau_{1}^{(i)} \leq t_{1}^{(j)} < t_{2}^{(j)} \leq \tau_{2}^{(i)}$.} \\ 
& \\
\hline
\multirow{2}{4em}{\centering \hyperref[C2]{C2}} & \multirow{2}{24em}{\centering $t_{1}^{(j)} \leq \tau_{1}^{(i)} < \tau_{2}^{(i)} \leq t_{2}^{(j)}$.} \\ 
& \\
\hline
\multirow{3}{4em}{\centering \hyperref[C3]{C3}} & \multirow{3}{24em}{\centering $t_{1}^{(j)} \leq \tau_{1}^{(i)} < t_{2}^{(j)} \leq \tau_{2}^{(i)}$ and $(t_{2}^{(j)} - \tau_{1}^{(i)}) > \rho_{\text{S}}(\tau_{2}^{(i)} - \tau_{1}^{(i)})$ or,
 $\tau_{1}^{(i)} \leq t_{1}^{(j)} < \tau_{2}^{(i)} \leq t_{2}^{(j)}$ and $(\tau_{2}^{(i)} - t_{1}^{(j)}) > \rho_{\text{S}}(\tau_{2}^{(i)} - \tau_{1}^{(i)})$.} \\ 
& \\
& \\
\hline
\hline
\multirow{2}{4em}{\centering \hyperref[D1]{D1}} & \multirow{2}{24em}{\centering $\tau_{2}^{(i)} \leq t_{2}^{(j)} < t_{1}^{(j+1)} \leq \tau_{1}^{(i+1)}$. } \\ 
& \\
\hline
\multirow{2}{4em}{\centering \hyperref[D2]{D2}} & \multirow{2}{24em}{\centering $\tau_{1}^{(i)} \leq t_{2}^{(j)} \leq \tau_{2}^{(i)}$ and $\tau_{1}^{(i+1)} \leq t_{1}^{(j+1)} \leq \tau_{2}^{(i+1)}$. } \\ 
& \\
\hline
\multirow{3}{4em}{\centering \hyperref[D3]{D3}} & \multirow{3}{34em}{\centering  
   $\tau_{1}^{(i)} \leq t_{2}^{(j)} \leq \tau_{2}^{(i)}$ and $\tau_{2}^{(i)} \leq t_{1}^{(j+1)} \leq \tau_{1}^{(i+1)}$ and $(t_{1}^{(j+1)} - \tau_{2}^{(i)}) > \rho_{\text{NS}}(t_{1}^{(j+1)} - t_{2}^{(j)})$ or, 
  $\tau_{2}^{(i)} \leq t_{2}^{(j)} \leq \tau_{1}^{(i+1)}$ and $\tau_{1}^{(i+1)} \leq t_{1}^{(j+1)} \leq \tau_{2}^{(j+1)}$ and $(\tau_{1}^{(i+1)} - t_{2}^{(j)}) > \rho_{\text{NS}}(t_{1}^{(j+1)} - t_{2}^{(j)})$.} \\ 
& \\
& \\
\hline
\end{tabular}
\caption{\label{tab:Case_i_a} Implementation of conditions used to classify seizure (\hyperref[C1]{C1}-\hyperref[C3]{C3}) and non-seizure (\hyperref[D1]{D1}-\hyperref[D3]{D3}) time intervals.}
\end{table}

\begin{table}[t]
\centering
\def\arraystretch{1.5}
\begin{tabular}{|c|c|c|c|c|c|c|c|}
\hline
Parameter setting & $\tau_{\text{NS}}$ & $\tau_{\text{S}}$ & $\tau_{w}$ & $\delta_{\alpha \beta}$ & $\alpha$ & $\beta$ & $\rho_{\text{S}}$ and $\rho_{\text{NS}}$ \\
\hline
$P_{1}$ & $0.7$ & $1$ & $0.5$ & $[0.01,0.06,0.01]$ & $[\delta_{\alpha \beta}+0.01,0.1, 0.005]$ & $\beta=\alpha-\delta_{\alpha \beta}$ & $[0,1,0.025]$  \\
\hline
$P_{2}$ & $2$ & $1$ & $1$ & " & " & " & "  \\
\hline
$P_{3}$ & $3$ & $2$ & $1$ & " & " & " & "  \\
\hline
\end{tabular}
\caption{\label{tab:AUC_params} 
Parameter settings used when evaluating the performance of the CT detection algorithm described in part \hyperref[M2]{M2} of the Methods section. Values specified in square brackets correspond to [lower value, upper value, difference between each value].
}
\end{table}


\subsection*{M4: Mitigating the occurrence of misleading points on ROC curves \label{M4}}

\begin{figure}[t]
    \centering
    \includegraphics[width=\textwidth]{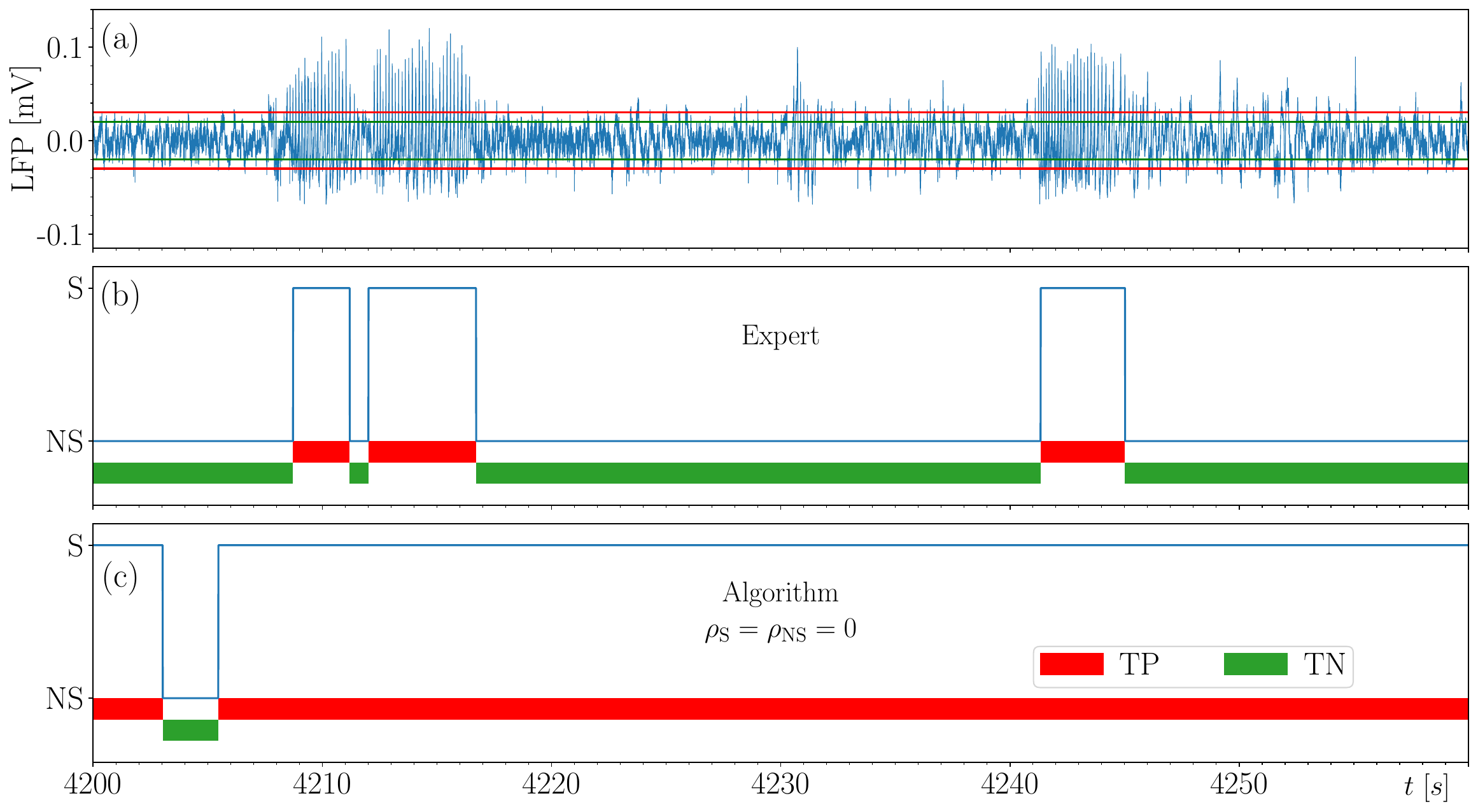}
    \caption{Example of misleading TP and TN classifications. (a) shows a portion of voltage recordings. (b) and (c) show binary sequence representations of the non-seizure and seizure time intervals in (a) according to expert annotations and our algorithm. The boxes beneath these sequences are coloured green for TN and red for TP. Algorithm parameters chosen as $\alpha=0.03$, $\beta=0.02$, $\tau_{\text{NS}}=2$, $\tau_{\text{S}}=1$, $\tau_{w}=1$, with $\rho_{\text{S}} = \rho_{\text{NS}}=0$ used in the classification of the seizure and non-seizure intervals.}
    \label{fig:MisClass_ROC_example1}
\end{figure}

In this subsection we (i) show that certain parameter settings of our algorithm can lead to misleading points on an ROC curve and (ii) present a method to prevent these misleading points from appearing. 
\\
\textbf{How misleading points occur: } 
Misleading points appear when the algorithm detects very few seizure and non-seizure time intervals in comparison to the expert and these time intervals satisfy the TP (\hyperref[C1]{C1}-\hyperref[C3]{C3}) and TN (\hyperref[D1]{D1}-\hyperref[D3]{D3}) conditions. 
As a result, the corresponding point on the ROC curve incorrectly indicates that the algorithm performed well.
This typically occurs when the upper voltage threshold, $\alpha$, is set too small or too large, causing the algorithm to detect one or few seizure time intervals whose total duration is disproportionally long or short in comparison to the time series the algorithm is applied to. For example, it can happen that a given $I_{\text{S}}^{(A)}$ contains almost all the $I_{\text{S}}^{(E)}$. 
We illustrate an example of such a scenario in Fig.~\ref{fig:MisClass_ROC_example1} when the algorithm is used to detect CTs between NS and S states in the T2M voltage recordings, a portion of which is shown in (a). 
In this example the algorithm parameters are chosen as $\alpha = 0.03$, $\beta = 0.02$, $\tau_{\text{NS}}=2$, $\tau_{\text{S}}=1$, $\tau_{w}=1$, and the time intervals are classified with overlap parameters set as $\rho_{\text{S}} = \rho_{\text{NS}}=0$. 
In (b) and (c) we show a binary sequence representation of the seizure and non-seizure time intervals according to expert annotations and the algorithm. Beneath each sequence we use the same coloured box convention in Fig.~\ref{fig:rhoS_effect_ROC_example1} to indicate how the corresponding seizure and non-seizure time intervals are classified. 
In this example, the expert identified 154 seizure time intervals, i.e., $N_{\text{S}}^{(E)}=154$, while the algorithm identified only $4$, i.e., $N_{\text{S}}^{(A)}=4$. However, these $4$ $I_{\text{S}}^{(A)}$ contain most of the $I_{\text{S}}^{(E)}$ and $I_{\text{NS}}^{(E)}$. 
Without adapting the classification procedure this results in (FPR, TPR) = (0,1), indicating that the algorithm achieves ideal performance. However, this is not the case in reality. 
\\
\textbf{Mitigating the occurrence of misleading points: } 
To combat the issue described above, we construct our ROC curve using (FPR, TPR) points  whose corresponding time intervals satisfy the following condition:
\begin{itemize}
    \item $\sum_{j=1}^{N_{\text{S}}^{(A)}} \left|\left(I_{\text{S}}^{(A)}\right)_{j}\right| \leq h_{u} \sum_{i=1}^{N_{\text{S}}^{(E)}} \left|\left(I_{\text{S}}^{(E)}\right)_{i}\right|$ for a suitably chosen $h_{u} > 0$. 
\end{itemize}
In other words, the sum of residence times in the S state according to the algorithm must be less than or equal to a multiple, $h_{u}$, of the sum of residence times in the S state according to the expert. 
It was found empirically that by setting $h_{u}=2$, this additional step to our classification procedure removes the misleading points on a given ROC curve.

\begin{figure}[t]
    \centering
    \includegraphics[width=\textwidth]{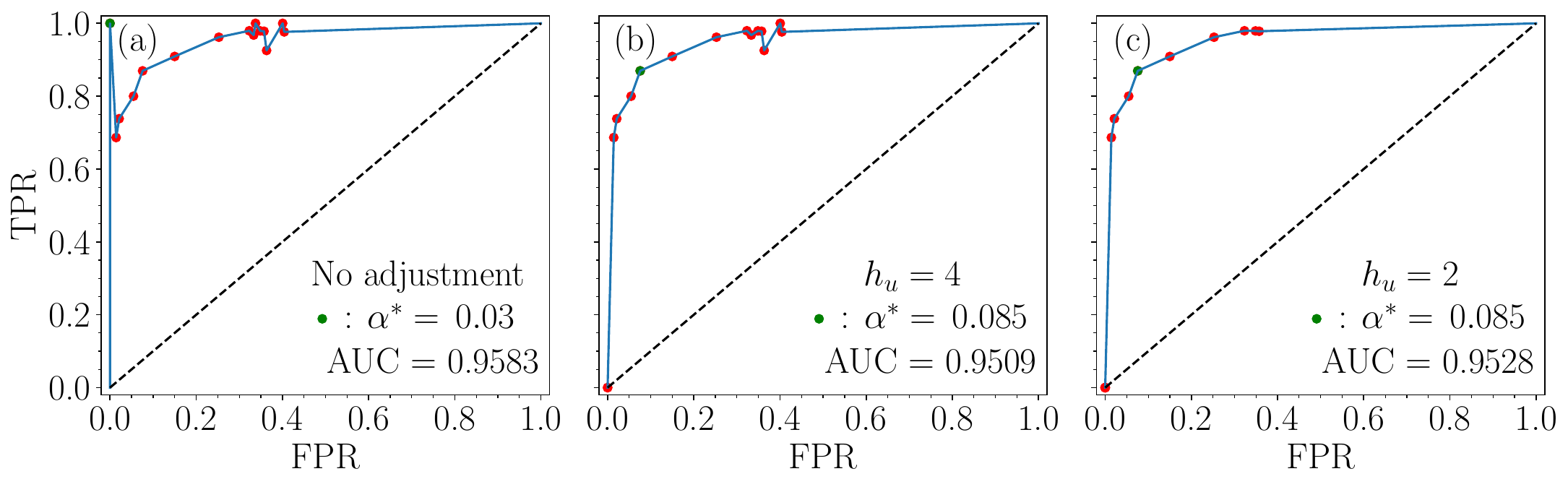}
    \caption{Effect of $h_{u}$ on removing misleading points from an ROC curve. (a) ROC curve without applying the $h_{u}$ adjustment. (b) and (c) ROC curves after applying the $h_{u}$ adjustment with $h_{u}=4$ and $h_{u}=2$, respectively. All ROC curves are generated using the same classifications of seizure and non-seizure time intervals in Fig.~\ref{fig:h2_best_rhoS_ROC_example1} for parameter setting $P_{2}$ with $\delta_{\alpha \beta} = 0.02$. In each panel, (FPR, TPR) points plotted in red correspond to different values of $\alpha$, the green point corresponds to $\alpha^{*}$ (the optimal $\alpha$), the lower right corner specifies $\alpha^{*}$ and the AUC.}
    \label{fig:MisClass_hoff_h30_h4_h2_ROC_example2}
\end{figure}

The result of applying this additional step for different choices of $h_{u}$ is illustrated in Fig.~\ref{fig:MisClass_hoff_h30_h4_h2_ROC_example2}. More specifically, we plot ROC curves based on applying our algorithm to recording session T2M with parameters setting $P_{2}$ and $\delta_{\alpha \beta}=0.02$ and classifying the resulting seizure and non-seizure time intervals with $\rho_{\text{S}} = \rho_{\text{NS}}=0$. 
Figure~\ref{fig:MisClass_hoff_h30_h4_h2_ROC_example2}~(a) shows the ROC curve obtained when no steps are taken to mitigate the occurrence of misleading points, the (FPR, TPR) point obtained from using $\alpha = 0.03$ incorrectly indicates the algorithm achieves ideal performance. 
Figure~\ref{fig:MisClass_hoff_h30_h4_h2_ROC_example2}~(b) shows the ROC curve obtained when setting $h_{u}=4$, the misleading points corresponding to $\alpha = 0.03$ and $0.035$ are removed. 
Figure~\ref{fig:MisClass_hoff_h30_h4_h2_ROC_example2}~(c) shows that by setting $h_{u}=2$, the misleading points corresponding to $\alpha \in \left[0.03,0.055\right]$ are removed. This also improves the overall smoothness of the ROC curve.
What is common across (a)-(c) is that the AUC is not highly sensitive to this method for $h \geq 2$. This is a desirable result as the sole motivation behind using this additional constraint was to mitigate the influence of these misleading points.

\section*{Acknowledgements}

This publication has emanated from research conducted with the financial support of Taighde \'{E}ireann – Research Ireland under grant number [19/FFP/6782].

\section*{Author contributions statement}

A.F. and S.W. conceived the experiments, A.F. conducted the experiments, A.F. and S.W. analysed the results, C.Mc~C., F.D., and V.C. provided the voltage recordings, C.Mc~C., K.L., and G.L. provided valuable discussions. A.F., C.Mc~C., K.L., G.L., and S.W. reviewed the manuscript. 

\section*{Additional information}

The authors have no competing interests to disclose.

\renewcommand{\thefigure}{S-\arabic{figure}}

\setcounter{figure}{0}

\newpage
\section*{Supplementary information\label{si:sec}}

\subsection*{S1: Wider study on how \texorpdfstring{$\tau_{\mathrm{NS}}$}{TEXT} and \texorpdfstring{$\tau_{\mathrm{S}}$}{TEXT} influence the algorithm's performance\label{S1}}

\begin{figure}
    \centering
    \includegraphics[width=.87\textwidth]{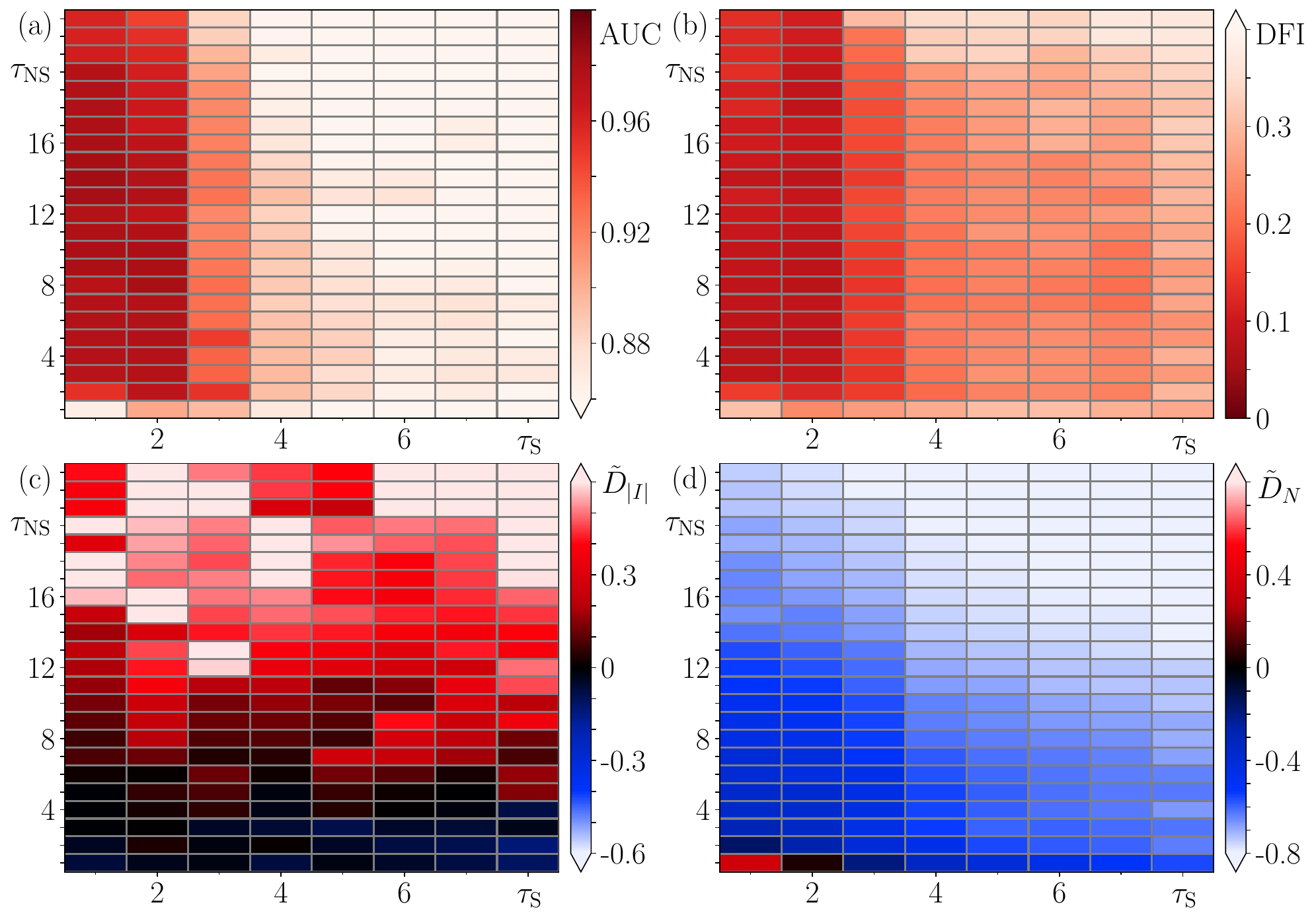}
    \caption{Extension to Fig.~\ref{fig:h2_best_rhoS_ROC_example1} in \hyperref[seref:MainText]{[1]}; illustrating how the values of $\tau_{\text{NS}}$ and $\tau_{\text{S}}$ influence the algorithm's performance in terms of the following metrics (a) AUC, (b) DFI, (c) $\widetilde{D}_{|I|}$, and (d) $\widetilde{D}_{|N|}$.} 
    \label{fig:_wider_study_tauNS_tauS_}
\end{figure}

In Fig.~\ref{fig:_wider_study_tauNS_tauS_} we extend the results shown in Fig.~\ref{fig:h2_best_rhoS_ROC_example1} in the main text to show how different choices of $\tau_{\text{S}} \in \left[1,8\right]$ and $\tau_{\text{NS}} \in \left[1,23\right]$ influence the algorithm's performance in terms of (a) the AUC, (b) the $\text{DFI}(\alpha^{*})$, (c) the relative difference between total time spent in the S state, denoted by $\widetilde{D}_{|I|}$ (see Eq.~\eqref{eq:DI}), and (d) the relative difference between total number of seizure time intervals, denoted by $\widetilde{D}_{N}$ (see Eq.~\eqref{eq:DN}). 
Note, we compute $\widetilde{D}_{|I|}$ and $\widetilde{D}_{N}$ based on $\left|I_{\text{S}}^{(A)}\right|$ and $N_{\text{S}}^{(A)}$ which are computed using the values of $t_{1}$ and $t_{2}$ obtained when using the optimal $\alpha$ (i.e., $\alpha^{*}$) for recording session T2M for a given choice of $\tau_{\text{S}}$ and $\tau_{\text{NS}}$. 
To maintain the connection to Fig.~\ref{fig:h2_best_rhoS_ROC_example1} in \hyperref[seref:MainText]{[1]}, we keep $\delta_{\alpha \beta} = 0.01$, $\tau_{w}=1$, and $\rho_{\text{S}}=\rho_{\text{NS}}=0$. Note, the tilde in $\widetilde{D}_{|I|}$ and $\widetilde{D}_{N}$ is to indicate that these are relative distances.

Figures~\ref{fig:_wider_study_tauNS_tauS_}~(a) and (b) show that while $\text{AUC}>0.97$ and $\text{DFI}<0.1$ for many choices of $\tau_{\text{NS}}$ when $\tau_{\text{S}}=1,2$. In contrast Figs.~\ref{fig:_wider_study_tauNS_tauS_}~(c) and (d) show that for larger $\tau_{\text{NS}}$ values, the total time spent in the S state according to the algorithm is significantly greater than according to the expert, and correspondingly, the algorithm detects significantly less seizure time intervals than the expert. 
Thus, the results presented in Fig.~\ref{fig:_wider_study_tauNS_tauS_} further inform our choice of $\tau_\text{S}=2$ and $\tau_\text{NS}=3$ to be optimal, reflecting a balance between performance and detecting a similar amount of seizure time intervals to the expert and that these intervals are of a similar lengths.
\begin{align}
    \widetilde{D}_{|I|} &= \dfrac{ \sum_{j=1}^{N_{\text{S}}^{(A)}} \left|\left(I_{\text{S}}^{(A)}\right)_{j}\right| - \sum_{i=1}^{N_{\text{S}}^{(E)}} \left|\left(I_{\text{S}}^{(E)}\right)_{i}\right| }{\sum_{i=1}^{N_{\text{S}}^{(E)}} \left|\left(I_{\text{S}}^{(E)}\right)_{i}\right|}.\label{eq:DI}\\
    \widetilde{D}_{N} &= \frac{ N_{\text{S}}^{(A)} - N_{\text{S}}^{(E)} }{N_{\text{S}}^{(E)}}.\label{eq:DN}
\end{align}


\subsection*{S2: Analysis of quantities used to generate ROC curves in Fig.~\ref{fig:h2_best_rhoS_ROC_example1} in \hyperref[seref:MainText]{[1]}\label{S2}}

\begin{figure}[t]
    \centering
    \includegraphics[width=.95\textwidth]{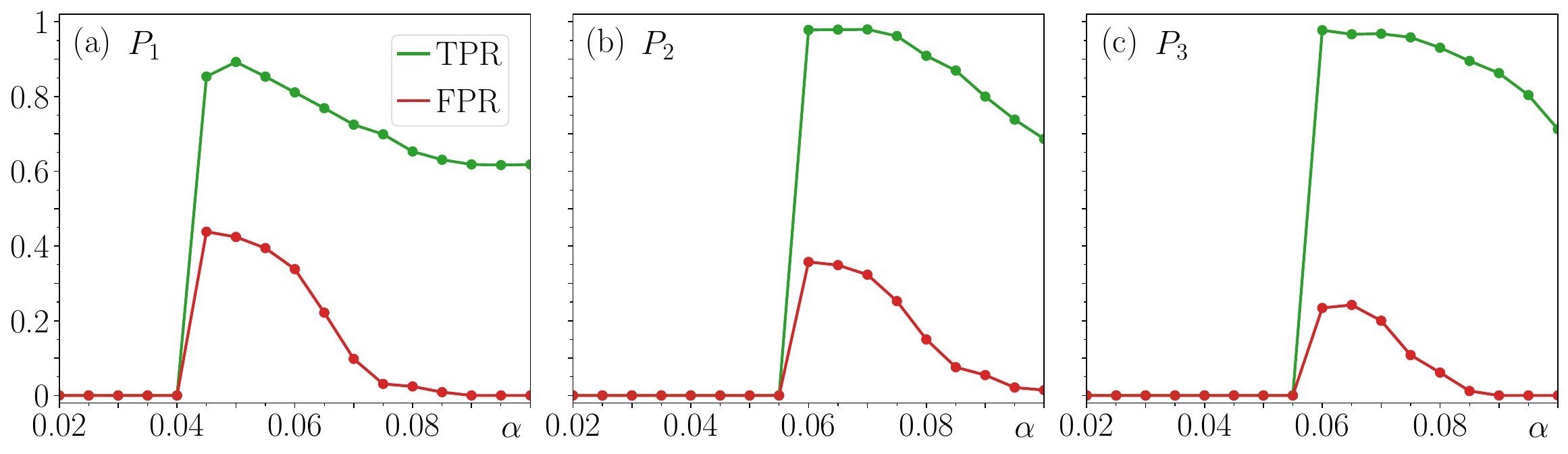}
    \caption{$\alpha$ vs. TPR (in green) and FPR (in red) for parameter setting $P_{1}$ in (a), $P_{2}$ in (b), and $P_{3}$ in (c).}
    \label{fig:P1P2P3_TPR_FPR_alpha_}
\end{figure}

\begin{figure}[t]
    \centering
    \includegraphics[width=.95\textwidth]{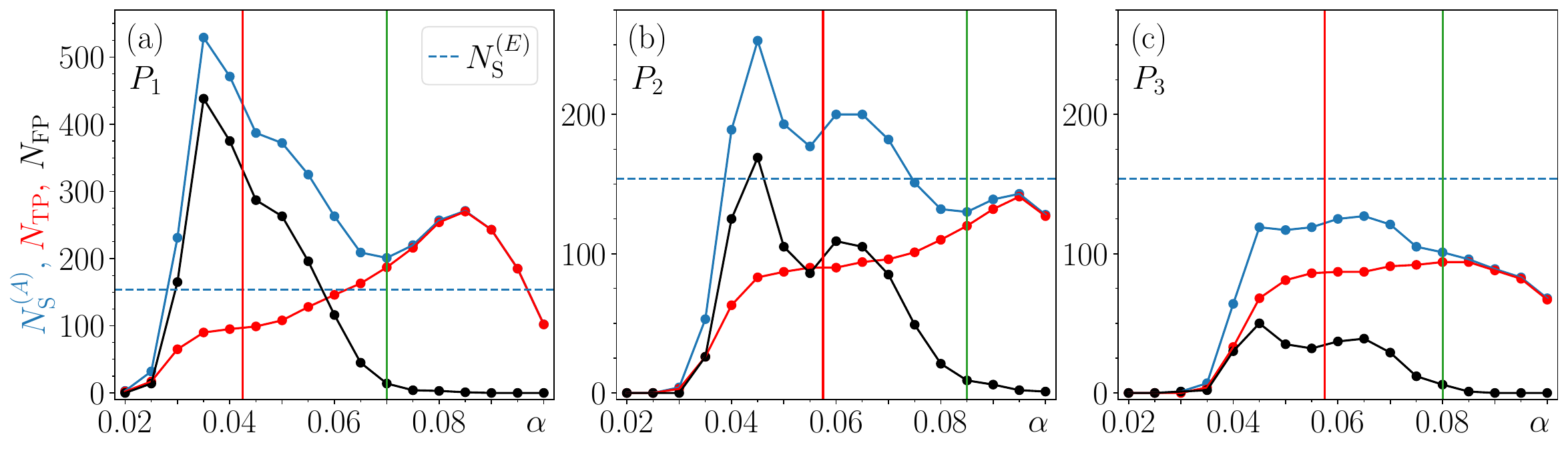}
    \caption{$\alpha$ vs. $N_{\text{S}}^{(A)}$ (in blue), $N_{\text{TP}}$ (in red), and $N_{\text{FP}}$ (in black) for parameter setting $P_{1}$ in (a), $P_{2}$ in (b), and $P_{3}$ in (c) and $\delta_{\alpha \beta} = 0.01$. The green vertical lines corresponds to the optimal $\alpha$. Dashed horizontal blue lines indicate $N_{\text{S}}^{(E)}$ for the T2M time series. Data points to the left of the red vertical line do not contribute to the corresponding ROC curves in Fig.~\ref{fig:h2_best_rhoS_ROC_example1} in \hyperref[seref:MainText]{[1]}.}
    \label{fig:S_alpha_vs_Ns_P123_deltaAB_0.01_T2M_}
\end{figure}

\begin{figure}[h!]
    \centering
    \includegraphics[width=.95\textwidth]{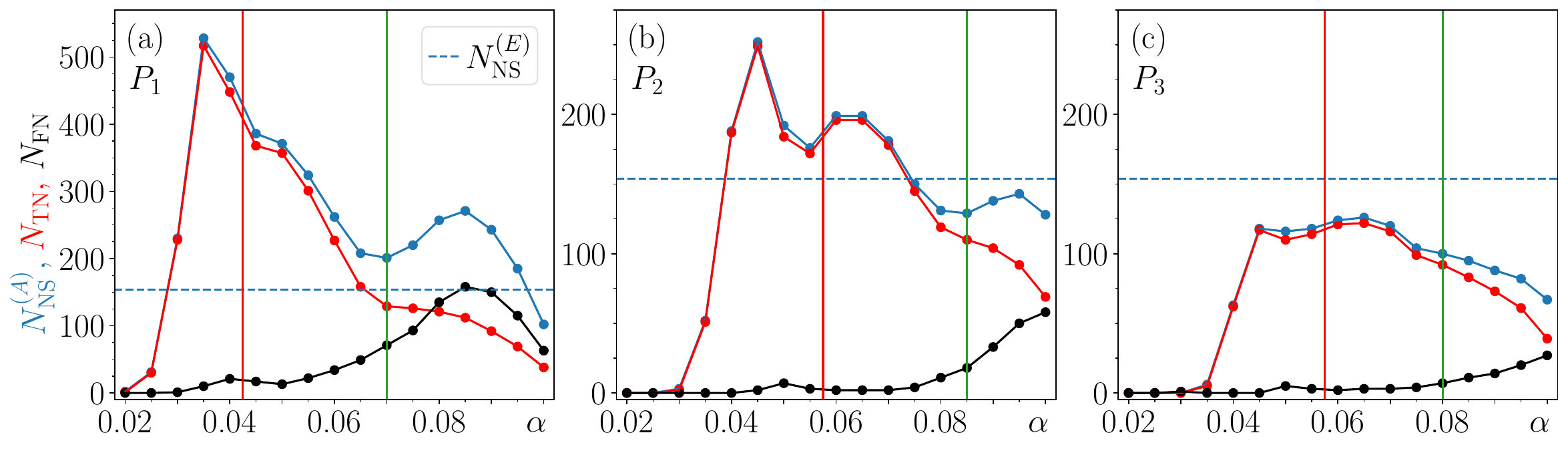}
    \caption{$\alpha$ vs. $N_{\text{NS}}^{(A)}$ (in blue), $N_{\text{TN}}$ (in red), and $N_{\text{FN}}$ (in black) for parameter setting $P_{1}$ in (a), $P_{2}$ in (b), and $P_{3}$ in (c) and $\delta_{\alpha \beta} = 0.01$. The green vertical lines corresponds to the optimal $\alpha$. Dashed horizontal blue lines indicate $N_{\text{NS}}^{(E)}$ for the T2M time series. Data points to the left of the vertical red line do not contribute to the corresponding ROC curves in Fig.~\ref{fig:h2_best_rhoS_ROC_example1} in \hyperref[seref:MainText]{[1]}.}
    \label{fig:NS_alpha_vs_Ns_P123_deltaAB_0.01_T2M_}
\end{figure}

Figures~\ref{fig:P1P2P3_TPR_FPR_alpha_}-\ref{fig:NS_alpha_vs_Ns_P123_deltaAB_0.01_T2M_} provide a more detailed breakdown on the different quantities that are used to construct the ROC curves in Fig.~\ref{fig:h2_best_rhoS_ROC_example1} in \hyperref[seref:MainText]{[1]}. 
Specifically, we show how these quantities vary with respect to $\alpha$ in parameter settings $P_1$, $P_2$, and $P_3$ in the respective panels (a), (b), and (c) in Figs.~\ref{fig:P1P2P3_TPR_FPR_alpha_}-\ref{fig:NS_alpha_vs_Ns_P123_deltaAB_0.01_T2M_}.
Figure~\ref{fig:P1P2P3_TPR_FPR_alpha_} shows how the TPR and FPR vary with respect to $\alpha$.
Figure~\ref{fig:S_alpha_vs_Ns_P123_deltaAB_0.01_T2M_} shows how $N_{\text{S}}^{(\text{A})}$, $N_{\text{TP}}$, and $N_{\text{FP}}$ vary with respect to $\alpha$ and how these quantities compare to $N_{\text{S}}^{(\text{E})}$.
Similarly, Figure~\ref{fig:NS_alpha_vs_Ns_P123_deltaAB_0.01_T2M_} shows how $N_{\text{NS}}^{(\text{A})}$, $N_{\text{TN}}$, and $N_{\text{FN}}$ vary with respect to $\alpha$ and how these quantities compare to $N_{\text{NS}}^{(\text{E})}$.
Further, in Figs.~\ref{fig:S_alpha_vs_Ns_P123_deltaAB_0.01_T2M_}-\ref{fig:NS_alpha_vs_Ns_P123_deltaAB_0.01_T2M_}, the purpose of the red vertical line is to emphasise that all information to the left of the line is not used to construct ROC curves, and the purpose of the green vertical line is to indicate the quantities which correspond to the optimal choice of $\alpha$, denoted by $\alpha^{*}$.

Figures~\ref{fig:P1P2P3_TPR_FPR_alpha_}-\ref{fig:NS_alpha_vs_Ns_P123_deltaAB_0.01_T2M_} show there is a nonlinear relationship between $\alpha$ and each quantity mentioned above.
We find that while the algorithm achieves its best performance for $P_3$, the algorithm detects less seizure time intervals in comparison to using parameter sets $P_2$ and $P_1$. This is to be expected since $\tau_{\text{S}}$ is smaller for $P_2$ and $P_1$. 
Furthermore, the motivation to exclude points from the ROC curves (see part \hyperref[M4]{M4} of the Methods section in \hyperref[seref:MainText]{[1]}) becomes clearer from Figs.~\ref{fig:S_alpha_vs_Ns_P123_deltaAB_0.01_T2M_} and \ref{fig:NS_alpha_vs_Ns_P123_deltaAB_0.01_T2M_} as these points correspond to cases where the algorithm detects a much larger number of seizure and non-seizure time intervals than the expert. 



\subsection*{S3: Assessing agreement beyond ROC curves \label{S3}}

In this subsection, we examine the agreement between our algorithm and the expert's annotations in ways beyond that captured by ROC curves, specifically the ROC curves shown in Fig.~\ref{fig:h2_best_rhoS_ROC_example1} in \hyperref[seref:MainText]{[1]}. 

\subsubsection*{S3 (a): Comparing $t_{1}$ and $\tau_{1}$}
\begin{figure}[t]
    \centering
    \includegraphics[width=.99\textwidth]{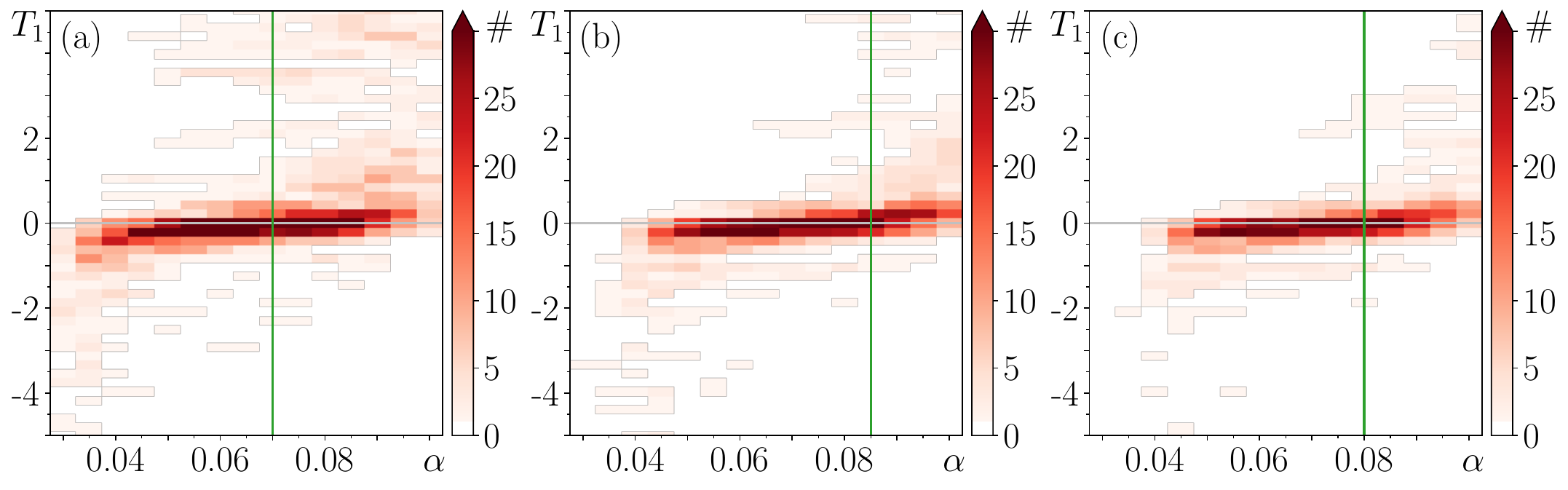}
    \caption{Heatmaps which show the number of TPs that take a particular $T_{1}=\tau_{1}-t_{1}$ value for a given $\alpha$ for parameter setting $P_{1}$ in (a), $P_{2}$ in (b), and $P_{3}$ in (c) and $\delta_{\alpha \beta} = 0.01$. The green vertical line in each panel corresponds to the optimal $\alpha$ for each parameter setting.}
    \label{fig:t1_vs_tau1_comparison_vs_alpha_P1_P2_P3_}
\end{figure}

In Fig.~\ref{fig:t1_vs_tau1_comparison_vs_alpha_P1_P2_P3_} we compare the times that our algorithm detects CTs from the NS to S state, $t_{1}$, to the expert annotation of seizure onset times, $\tau_{1}$, for $I_{\text{S}}^{(\text{A})}$ that were classified as TPs when using parameter setting $P_1$ in (a), $P_2$ in (b), and $P_3$ in (c). 
For each TP, we compute the time difference $T_{1}=\tau_{1}-t_{1}$ and show, for each value of $\alpha$, the number of TPs associated with each $T_{1} \in \left[ -5, 5 \right]$. Darker shades of red indicate a larger number of TPs.

Across Figs.~\ref{fig:t1_vs_tau1_comparison_vs_alpha_P1_P2_P3_}~(a)-(c) we observe the following common trends: 
$T_{1}<0$ for most TPs when $\alpha \lesssim 0.055$, meaning that our algorithm typically detects CTs from the NS to S state at times before the expert's annotation of the corresponding seizure onset. 
$T_{1}\approx0$ for most TPs when $0.055 \lesssim \alpha \lesssim 0.085$, thus showing close agreement between our algorithm and the expert annotations, with $\tau_{1}$ and $t_{1}$ differing by $0.5\text{s}$ in most cases.
$T_{1} > 0$ for most TPs when $\alpha \gtrsim 0.085$, indicating that our algorithm typically detects most CTs from the NS to S state at times after the expert's annotation of the corresponding seizure onset. 

\subsubsection*{S3 (b): Comparing residence times}
\begin{figure}[t]
    \centering
    \includegraphics[width=.8\textwidth]{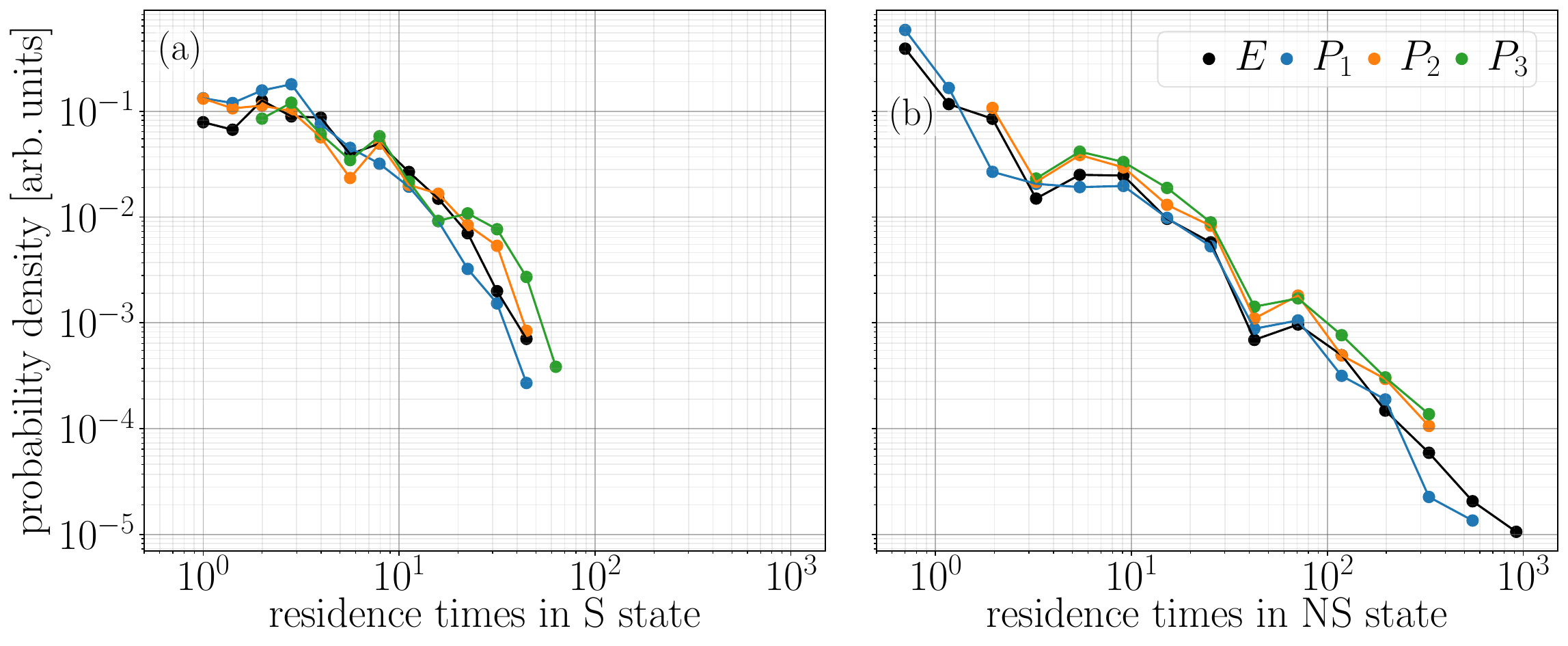}
    \caption{Probability density of residence times in (a) the S state and (b) the NS state based on seizure and non-seizure time intervals according to the expert (denoted by E) and the algorithm in the parameter settings $P_{1}$, $P_{2}$, and $P_{3}$ for $\delta_{\alpha \beta} = 0.01$, computed for the optimal $\alpha$ in parameter setting.}
    \label{fig:ResTimes_CianAlg}
\end{figure}

Figure~\ref{fig:ResTimes_CianAlg} compares the probability densities of residence times in the (a) S and (b) NS states obtained from our algorithm with those derived from the expert annotations. 
Specifically, we compare the lengths of the different $\left(I_{\text{S}}^{(E)}\right)_{i}$ and $\left(I_{\text{NS}}^{(E)}\right)_{i}$ to $\left(I_{\text{S}}^{(A)}\right)_{j}$ and $\left(I_{\text{NS}}^{(A)}\right)_{j}$ for the T2M time series studied in \hyperref[seref:MainText]{[1]}.
The residence times derived from the experts annotations are shown in black and the algorithm-derived residence times are computed from the $t{1}$ and $t_{2}$ values obtained when using the optimal $\alpha$ values, $\alpha^{*}$, for parameter settings $P_{1}$ (in blue), $P_{2}$ (in orange), and $P_{3}$ (in green).

Overall, the residence times obtained from our algorithm closely match those obtained from the expert annotations for both the S and NS states across all three parameter settings.

\subsection*{S4: Analysis of PV curves (complimentary to analysis of ROC curves)\label{S4}}

\begin{figure}[t]
    \centering
    \includegraphics[width=.95\textwidth]{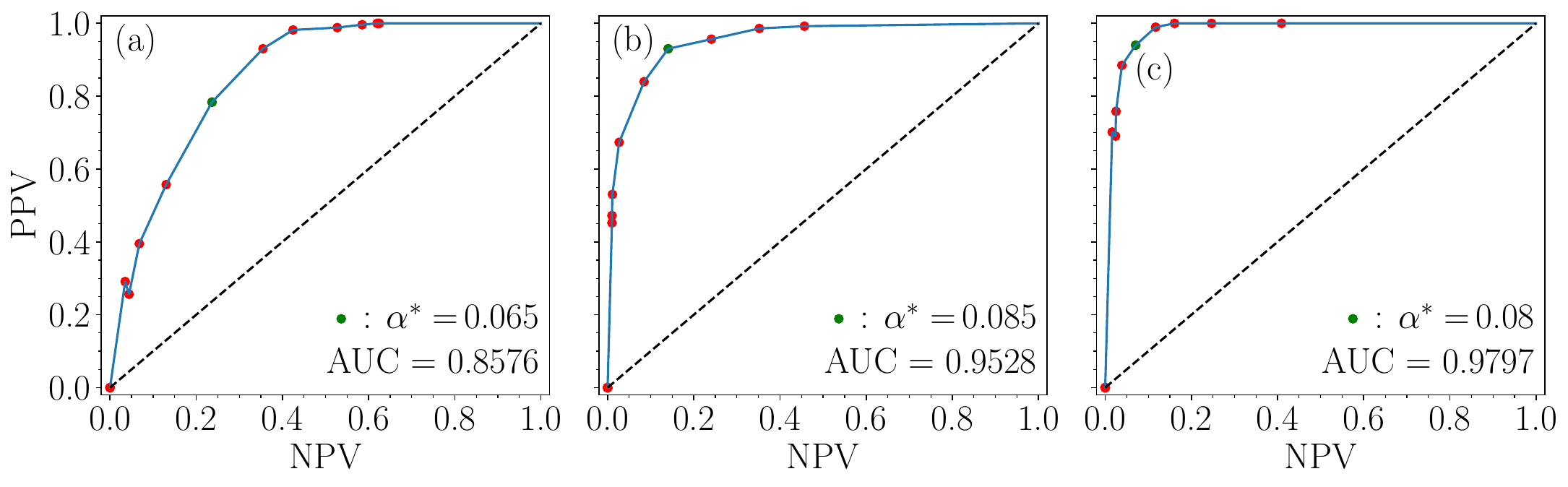}
    \caption{PV curves obtained when using parameter setting $P_1$ in (a), $P_2$ in (b), and $P_3$ in (c). In each panel, (NPV, PPV) points plotted in red correspond to different values of $\alpha$, the green point corresponds to $\alpha^{*}$ (the optimal $\alpha$), the lower right corner specifies $\alpha^{*}$ and the AUC.}
    \label{fig:h2_best_rhoS_PV_example1}
\end{figure}

Taking further inspiration from Temko \textit{et al.}~\hyperref[seref:Temko]{[2]} and Mathieson \textit{et al.}~\hyperref[seref:Mathieson]{[3]}, we also consider the following metrics to examine the agreement between our algorithm and the experts annotations:
\\
\textbf{Positive predictive value (PPV):} $N_{\text{TP}} / N_{\text{S}}^{(A)} \in \left[ 0, 1 \right]$, also referred to as the `seizure detection rate', defines the proportion of correctly classified seizure time intervals. 
\\
\textbf{Negative predictive value (NPV):} $N_{\text{FP}} / N_{\text{NS}}^{(A)} \in \left[ 0, 1 \right]$, the proportion of incorrectly classified non-seizure time intervals. 
\\
Note, the NPV is a modified version of the `false detections per hour' (FD/h) metric considered by Temko \textit{et al.} and Mathieson \textit{et al.} which is used to quantify the number of false positives that occur per hour in a given time series.
The FD/h metric is an insightful clinical metric when the time series that are analysed are several hours long. 
Since the time series we analyse are much shorter, and are at most 1-2 hours long, we use the NPV to quantify similar behaviour. 
\\
\\
We use the PPV and NPV metrics to construct what we call `predictive value' (PV) curves using the same procedure as ROC curves specified in \hyperref[seref:MainText]{[1]} for increasing values of NPV. 
Similar to the TPR and FPR, the closer the PPV is to $1$ and the NPV is to $0$ the better the performance of our algorithm. 
In Fig.~\ref{fig:h2_best_rhoS_PV_example1} we show PV curves for the same parameter settings as in Fig.~\ref{fig:h2_best_rhoS_ROC_example1} in \hyperref[seref:MainText]{[1]}. Each PV curve provides similar insight to its ROC counterpart. 

\section*{References}

\textbf{[1]. \label{seref:MainText}} Flynn, A. \textit{et al.} ``Detecting seizure onset and offset times using human intelligence: A critical transitions based approach”.\\
\textbf{[2]. \label{seref:Temko}} Temko, A. \textit{et al.} ``Performance assessment for EEG-based neonatal seizure detectors.'' \textit{Clin. Neurophysiol.} \textbf{122}, 474–482 (2011).\\
\textbf{[3]. \label{seref:Mathieson}} Mathieson, S. R. \textit{et al.} ``Validation of an automated seizure detection algorithm for term neonates.'' \textit{Clin. Neurophysiol.} \textbf{127}, 156–168 (2016).

\end{document}